\documentclass[onefignum,onetabnum]{siamart251104}

\usepackage{microtype}

\usepackage{xcolor}

\usepackage{amsmath,amssymb} % AMS系，ほぼ必須
\usepackage{mathrsfs} % 花文字フォント (\mathscr)
\usepackage{url} % url 用フォント (\url)
\usepackage{latexsym} % マニアックな数学記号用
\usepackage{bm} % 斜体太字フォント (\bm)

\usepackage{mathtools}
\usepackage{booktabs,multirow}
\usepackage{array}
\newcolumntype{P}[1]{>{\centering\arraybackslash}p{#1}}
\usepackage{diagbox}
\usepackage{threeparttable}
\usepackage{siunitx} % for scientific notation
\usepackage{tabularx}
\newcolumntype{C}{>{\centering\arraybackslash}X}

\usepackage{algorithm}
\usepackage[noend]{algpseudocode}

\algnewcommand{\Break}{\textbf{break}}
\algnewcommand{\algorithmicgoto}{\textbf{go to} Line}
\algnewcommand{\Goto}[1]{\algorithmicgoto~\ref{#1}}
\algblockdefx{MRepeat}{EndRepeat}{\textbf{repeat}}{}
\algnotext{EndRepeat}
\usepackage[sort&compress,numbers]{natbib}
\usepackage{thmtools,thm-restate}

\usepackage{hyperref}
\usepackage{cleveref}
\expandafter\def\csname ver@etex.sty\endcsname{3000/12/31} % autonum のバグ対策
\usepackage{autonum} % 参照されたときのみ式番号を表示 (cleveref と共存可)
\makeatletter
\autonum@generatePatchedReferenceCSL{Cref} % autonum + cleveref のバグ対策
\makeatother

\usepackage{multicol}
\usepackage{adjustbox}
\usepackage{placeins} % FloatBarrier to control float placement
\usepackage{xparse} % 複雑なマクロ定義用
\usepackage{pbox} % pbox コマンド用
\usepackage{xspace} % newcommand マクロ直後に適切にスペースを入れる

\usepackage[subrefformat=parens]{subcaption}
\usepackage{enumitem} % enumerate package より詳細な設定
\DeclareMathOperator{\diag}{diag}

\DeclareMathOperator{\tr}{tr}

\DeclareMathOperator*{\tsum}{\textstyle\sum}

\DeclareMathOperator{\EE}{\mathbb{E}}

\DeclarePairedDelimiter\floor{\lfloor}{\rfloor}
\DeclarePairedDelimiterX\inner[2]{\langle}{\rangle}{{#1},{#2}}
\DeclarePairedDelimiter\abs{\lvert}{\rvert}
\DeclarePairedDelimiter\norm{\lVert}{\lVert}
\DeclarePairedDelimiter\normop{\lVert}{\lVert_{\mathrm{op}}}
\DeclarePairedDelimiter\normFro{\lVert}{\lVert_{\mathrm{F}}}
\DeclarePairedDelimiter\set{\{}{\}}
\DeclarePairedDelimiter\prn{(}{)}
\DeclarePairedDelimiter\bra{[}{]}
\DeclarePairedDelimiterX\Set[2]{\{}{\}}{\mspace{2mu}{#1}\;\delimsize|\;{#2}\mspace{2mu}}
\DeclarePairedDelimiterX\Prn[2]{(}{)}{\mspace{2mu}{#1}\;\delimsize|\;{#2}\mspace{2mu}}
\DeclarePairedDelimiterX\Bra[2]{[}{]}{\mspace{2mu}{#1}\;\delimsize|\;{#2}\mspace{2mu}}

\newcommand{\N}{\mathbb N}

\newcommand{\R}{\mathbb R}

\newcommand{\0}{\mathbf 0}

\newcommand{\e}{\mathrm{e}}
\newcommand{\dd}{\mathrm{d}}
\renewcommand{\O}{\mathrm{O}}

\renewcommand{\epsilon}{\varepsilon}
\NewDocumentCommand{\exsub}{s m O{} m}{%
  \IfBooleanT{#1}{\EE_{#2}\nolimits\bra*{#4}}%
  \IfBooleanF{#1}{\EE_{#2}\nolimits\bra[#3]{#4}}%
}
\NewDocumentCommand{\ex}{s O{} m}{%
  \IfBooleanT{#1}{\EE\nolimits\bra*{#3}}%
  \IfBooleanF{#1}{\EE\nolimits\bra[#2]{#3}}%
}
\NewDocumentCommand{\cex}{s O{} m m}{%
  \IfBooleanT{#1}{\EE\nolimits\Bra*{#3}{#4}}%
  \IfBooleanF{#1}{\EE\nolimits\Bra[#2]{#3}{#4}}%
}
\usepackage{CJKutf8}

\newcommand{\mathInd}{\hphantom{{}={}}}
\newcommand{\by}[2][]{\text{\pbox[c]{\textwidth}{(by \pbox[t]{\textwidth}{\phantom{}#2)#1}}}}

\newcommand{\xinit}{x_{\mathrm{init}}}

\renewcommand{\SS}{\mathbb{S}}

\newcommand{\Proposed}{\textsf{Proposed}\xspace}
\newcommand{\BFGS}{\textsf{BFGS}\xspace}
\newcommand{\DFP}{\textsf{DFP}\xspace}
\newcommand{\GD}{\textsf{GD}\xspace}
\newcommand{\AGDR}{\textsf{AGDR}\xspace}

\usepackage{lipsum}
\usepackage{amsfonts}
\usepackage{graphicx}
\usepackage{epstopdf}
\ifpdf
  \DeclareGraphicsExtensions{.eps,.pdf,.png,.jpg}
\else
  \DeclareGraphicsExtensions{.eps}
\fi

\newsiamremark{remark}{Remark}
\newsiamthm{assumption}{Assumption}

\crefname{section}{Section}{Sections}
\Crefname{section}{Section}{Sections}
\crefname{subsection}{Section}{Sections}
\Crefname{subsection}{Section}{Sections}

\setlist[itemize]{
  topsep=0.1\baselineskip,
  itemsep=0\baselineskip,
  leftmargin=1.5em,
}

\newlist{enuminasm}{enumerate}{1} % asm 環境内で enumerate を使うときの環境
\setlist[enuminasm]{
  font=\upshape,
  label=(\alph*),
  ref=\theassumption(\alph*),
  topsep=0.2\baselineskip,
  itemsep=0.1\baselineskip,
  leftmargin=2em,
}
\crefalias{enuminasmi}{assumption}

\headers{Parameter-Free Accelerated Quasi-Newton}{N. Marumo}

\title{Parameter-Free Accelerated Quasi-Newton Method for Nonconvex Optimization\thanks{Submitted to the editors on December 10, 2025.
\funding{This work was funded by JSPS KAKENHI (24K23853) and JST CREST (JPMJCR24Q2).}}}

\author{Naoki Marumo\thanks{University of Tokyo, Tokyo, Japan 
  (\email{marumo@mist.i.u-tokyo.ac.jp}).}}

\usepackage{amsopn}

\ifpdf
\hypersetup{
  pdftitle={Parameter-Free Accelerated Quasi-Newton Method for Nonconvex Optimization},
  pdfauthor={N. Marumo}
}
\fi

\begin{document}

\maketitle

% REQUIRED
\begin{abstract}
We propose a quasi-Newton-type method for nonconvex optimization with Lipschitz continuous gradients and Hessians. The algorithm finds an $\varepsilon$-stationary point within $\mathrm{O}(d^{1/4} \varepsilon^{-13/8})$ function and gradient evaluations, where $d$ is the problem dimension. Our method is \emph{parameter-free} in the sense that it requires no prior knowledge of problem-dependent parameters such as Lipschitz constants or the optimal value. Moreover, it does not need the target accuracy $\varepsilon$ or the total number of iterations to be specified in advance. The result is achieved by combining several key ideas: momentum-based acceleration, quartic regularization for subproblems, and a scaled variant of the Powell-symmetric-Broyden (PSB) update.
\end{abstract}

% REQUIRED

% REQUIRED

\section{Introduction}
We consider the general nonconvex optimization problem
\begin{align}
  \min_{x \in \R^d} \ f(x),
  \label{eq:main_problem}
\end{align}
where $f \colon \R^d \to \R$ is twice differentiable and bounded below, i.e., $f^* \coloneqq \inf_{x \in \R^d} f(x) > -\infty$.
Throughout the paper, we assume that $f$ has Lipschitz continuous gradients and Hessians, as formalized below.
\begin{assumption}
  \label{asm:lipschitz_grad_hess}
  Let $L, M > 0$ be constants.
  \begin{enuminasm}
    \item 
    \label{asm:lipschitz_grad}
    $\norm*{\nabla f(x) - \nabla f(y)} \leq L \norm*{x - y}$ for all $x, y \in \R^d$.
    \item
    \label{asm:lipschitz_hess}
    $\normop{\nabla^2 f(x) - \nabla^2 f(y)} \leq M \norm*{x - y}$ for all $x, y \in \R^d$.
  \end{enuminasm}
\end{assumption}

First-order methods, which access $f$ only through function values and gradients, are fundamental tools in optimization and its applications.
The classical gradient descent method is known to find an $\epsilon$-stationary point (i.e., a point $x \in \R^d$ at which $\norm{\nabla f(x)} \leq \epsilon$) within $\O(\epsilon^{-2})$ gradient evaluations under the Lipschitz continuity of gradients (e.g., \citep[Section~1.2.3]{nesterov2004introductory}).
By additionally assuming Lipschitz continuous Hessians, many sophisticated first-order methods have been proposed to achieve improved complexity bounds of $\tilde{\O}(\epsilon^{-7/4})$ or $\O(\epsilon^{-7/4})$ \citep{carmon2017convex,jin2018accelerated,allen2018neon2,xu2017neon+,cutkosky2023optimal,li2022restarted,li2023restarted,marumo2024parameter,marumo2025universal}.

A significant recent advance is due to Jiang et al.~\citep{jiang2025improved}, who introduced a first-order algorithm based on a quasi-Newton method that achieves the dimension-dependent complexity $\O(d^{1/4} \epsilon^{-13/8})$.
This bound improves upon $\O(\epsilon^{-7/4})$ when $d = \O(\epsilon^{-1/2})$, which corresponds to medium-scale regimes.
However, the algorithm requires the knowledge of problem-dependent parameters such as the Lipschitz constants $L$ and $M$ and the optimal value $f^*$, and it also needs the total number of iterations to be specified in advance.
These quantities are often unavailable in practice, which limits the applicability of the method.

To address these issues, we propose a new quasi-Newton-type method with the following advantages.
\begin{itemize}
  \item 
  The algorithm finds an $\epsilon$-stationary point within $\O(d^{1/4} \epsilon^{-13/8})$ function and gradient evaluations, matching the best known dimension-dependent complexity.
  \item
  It is \emph{parameter-free} in the sense that it requires no prior knowledge of problem-dependent parameters such as $L$, $M$, or $f^*$.
  \item
  It does not require the target accuracy $\epsilon$ or the total number of iterations to be specified in advance.
  \item
  The algorithm evaluates a gradient only once per iteration.
\end{itemize}

To achieve these advantages, our method is built on two key ideas: a quartic-regularized subproblem and a scaled variant of the classical Powell-symmetric-Broyden (PSB) update~\citep{powell1970new}.
Quasi-Newton methods update the iterate $x_k$ and the Hessian approximation $B_k$, and both updates in \citep{jiang2025improved} involve parameters that depend on Lipschitz constants.
To eliminate the need to know the Lipschitz constants, we instead adopt a quartic-regularized subproblem for updating $x_k$ together with a scaled PSB update for $B_k$.
In addition, we incorporate an acceleration technique, inspired by momentum methods \citep{polyak1964methods,nesterov1983method}, to achieve the desired complexity bound.

We evaluate the numerical performance of the proposed method on benchmark problems and on an inverse problem for a nonlinear wave equation.
The results indicate that our algorithm is competitive with baseline methods.

\paragraph{Notation}
Let $\N \coloneqq \set{0, 1, 2, \ldots}$.
Let $\R^d$ be the $d$-dimensional Euclidean space equipped with the standard inner product $\inner*{\cdot}{\cdot}$ and the induced norm $\norm*{\cdot}$.
Let $\SS^d$ denote the set of $d \times d$ real symmetric matrices.
For $A, B \in \SS^d$, let $\inner*{A}{B} \coloneqq \tr (A^\top B)$.
We denote by $\normop*{\cdot}$ the operator (spectral) norm of a matrix and by $\normFro*{\cdot}$ its Frobenius norm.

\section{Related work}
This section reviews related work from several perspectives and highlights the differences from the present work.

\paragraph{Complexity bounds of first-order methods}
Carmon et al.~\citep{carmon2017convex} first achieved the complexity $\tilde{\O}(\epsilon^{-7/4})$ using a first-order method under \cref{asm:lipschitz_grad_hess}.\footnote{
  Agarwal et al.~\citep{agarwal2017finding} and Carmon et al.~\citep{carmon2018accelerated} also achieved the complexity $\tilde{\O}(\epsilon^{-7/4})$, but their algorithms rely on Hessian-vector products and are therefore not purely first-order.
}
Several subsequent algorithms with the same complexity have since been proposed \citep{jin2018accelerated,allen2018neon2,xu2017neon+,cutkosky2023optimal}.
Li and Lin~\citep{li2022restarted,li2023restarted} removed the logarithmic factor and obtained the complexity $\O(\epsilon^{-7/4})$, while parameter-free variants with the same complexity were also developed \citep{marumo2024parameter,marumo2025universal}.
All these results provide bounds that do not depend on the dimension $d$.
In contrast, some first-order methods~\citep{cartis2012oracle,grapiglia2022cubic,cartis2025global} have dimension-dependent complexity bounds of $\O(d \epsilon^{-3/2})$.
This complexity has been improved to $\O(\sqrt{d} \epsilon^{-3/2})$~\citep{doikov2025first}.
Jiang et al.~\citep{jiang2025improved} developed a first-order method with a complexity bound of $\O(d^{1/4}\epsilon^{-13/8})$.
In this paper, we propose a parameter-free first-order method whose complexity matches that of \citep{jiang2025improved}.
A detailed comparison with \citep{jiang2025improved} is provided in \cref{sec:comparison_with_jiang}.

\paragraph{Update rules of Hessian approximations}
Quasi-Newton methods generate a sequence of Hessian approximations $(B_k)_{k \in \N}$, and various update rules have been developed, including BFGS, DFP, and SR1 (see, e.g., \citep[Chapter~6]{nocedal2006numerical}).
BFGS and DFP are highly successful for strongly convex optimization, and their non-asymptotic superlinear convergence rates have been recently established~\citep{rodomanov2022rates,jin2025non,jin2023non}.
However, their update rules contain the denominator $\inner*{\nabla f(x_{k+1}) - \nabla f(x_k)}{x_{k+1} - x_k}$, which may vanish in nonconvex settings.
In contrast, the Powell-symmetric-Broyden (PSB) formula~\citep{powell1970new}, whose denominator is $\norm*{x_{k+1} - x_k}^2$, is more suitable for nonconvex problems.
Cartis et al.~\citep{cartis2025global} analyzed a cubic-regularized PSB method for nonconvex optimization and established a complexity bound of $\O(d \epsilon^{-3/2})$.
Our method also uses the PSB update, but incorporates a scaling factor.

\paragraph{Parameter-free algorithms}
Removing the dependence on unknown problem parameters, such as Lipschitz constants, has long been a central theme in algorithm design.
A common strategy is to maintain an estimate of a Lipschitz constant and increase it until a verifiable descent condition is satisfied, as in Armijo's backtracking~\citep{armijo1966minimization} and trust-region methods (see, e.g., \citep{conn2000trust}).
In our quasi-Newton scheme, however, the relevant descent inequality involves the true Hessian and is therefore inaccessible to the algorithm, thus making these approaches difficult to employ.
Many other strategies have been developed, including the adaptive Polyak step size~\citep{hazan2019revisiting}, adaptive gradient descent~\citep{malitsky2020adaptive}, AdaGrad-Norm~\citep{levy2017online,ward2020adagrad}, and AdaBB~\citep{zhou2025adabb}.
These methods are designed for gradient-descent-type updates or convex objectives, and extending them to nonconvex quasi-Newton schemes is nontrivial.
To bypass these limitations, we adopt a different approach: combining quartic regularization with a scaled PSB update.

\paragraph{Higher-order regularization}
Regularization for subproblems is often quadratic or cubic (e.g., \citep{ueda2010convergence,mishchenko2023regularized,nesterov2006cubic,doikov2024super}), yet higher-order regularization has also been explored.
Birgin et al.~\citep{birgin2017worst} analyzed $p$th-order methods using a $(p+1)$st-order regularized subproblem and established the complexity $\O(\epsilon^{-\frac{p+1}{p}})$.
Assuming $\nu$-H\"older continuous gradients, Cartis et al.~\citep{cartis2017worst} achieved the complexity $\O(\epsilon^{-\frac{1+\nu}{\nu}})$ by using $r$th-order regularization with $r > 1 + \nu$ within a first-order method.
These results were unified in \citep{cartis2019universal} to obtain general complexity bounds for higher-order methods.
Cartis and Zhu~\citep{cartis2025second} used quartic-regularized quadratic minimization in the context of third-order methods.
Our algorithm also employs quartic regularization, but the motivation is different: quartic regularization helps ensure sufficient decrease of the objective value when combined with a scaled PSB update, as discussed in \cref{rem:why_quartic}.

\section{Proposed algorithm}
\label{sec:proposed_algorithm}

The proposed algorithm is an accelerated quasi-Newton method, as summarized in \cref{alg:proposed}.
The algorithm has a double-loop structure.
In outer iteration $t$, we fix four parameters $\kappa, \sigma, \delta, \theta > 0$ and then repeat accelerated quasi-Newton updates for $K \coloneqq \floor*{\kappa}$ inner iterations, thereby generating the iterates $(x_k)_{k=0}^K$.

For $k \geq 1$, let $\bar x_k$ and $\bar g_k$ be the weighted averages defined as
\begin{align}
  \bar x_k
  &\coloneqq
  \frac{1}{k (k+1)} \prn*{ \sum_{i=0}^{k-1} (2i+1) x_i + k x_k },
  \label{eq:def_xbar}\\
  \bar g_k
  &\coloneqq
  \frac{1}{k (k+1)} \prn*{ \sum_{i=0}^{k-1} (2i+1) \nabla f(x_i) + k \nabla f(x_k) }.
  \label{eq:def_gbar}
\end{align}
As shown in \cref{sec:complexity_analysis}, the algorithm ensures that $\norm*{\nabla f(\bar x_K)} < \epsilon$ after sufficiently many iterations.
Although $\bar g_k$ does not explicitly appear in the algorithm, it serves as a proxy for $\nabla f(\bar x_k)$ in the analysis.
The weights are chosen so that
\begin{align}
  \label{eq:relation_gbar_g}
  \ \ 
  (k+2) \bar g_{k+1} - k \bar g_k
  &=
  \frac{1}{k+1} \prn*{ \sum_{i=0}^k (2i+1) \nabla f(x_i) + (k+1) \nabla f(x_{k+1}) }\\
  &\qquad
  - \frac{1}{k+1} \prn*{ \sum_{i=0}^{k-1} (2i+1) \nabla f(x_i) + k \nabla f(x_k) }\\
  &=
  \frac{1}{k+1} \prn[\big]{ (2k+1) \nabla f(x_k) + (k+1) \nabla f(x_{k+1}) - k \nabla f(x_k) }\!\!\!\!\\
  &=
  \nabla f(x_k) + \nabla f(x_{k+1}),
\end{align}
which is a key relation for proving sufficient decrease of the objective value.

\begin{algorithm}[t]
  \caption{Proposed accelerated quasi-Newton method\label{alg:proposed}}
  \begin{algorithmic}[1]
    \Require{%
      $\xinit \in \R^d$,
      $((\kappa_t, \sigma_t, \delta_t))_{t \in \N}$,
      where $\kappa_t > d^{1/5}$ and $\sigma_t, \delta_t > 0$
      for all $t \in \N$
    }
    \State{%
      $x_0 \gets \xinit$,
      $B_0 \gets O$
    }
    \For{$t=0,1,\dots$}
      \State{%
        $(\kappa, \sigma, \delta) \gets (\kappa_t, \sigma_t, \delta_t)$,
        $\theta \gets d / \kappa^5$,
        $K \gets \floor*{\kappa}$
      }
      \label{alg-line:set_parameters}
      \For{$k=0,1,\dots,K-1$}
        \label{alg-line:inner_loop_start}
        \State $m_k(s) \coloneqq \inner[\big]{\nabla f(x_k) + \frac{1}{k+1} \tsum_{i=0}^k (2i+1) \nabla f(x_i)}{s} + \frac{1}{2} \inner*{B_k s}{s} + \frac{\sigma}{4} \norm*{s}^4$
        \State $s_k \in \Set{s \in \R^d}{\norm*{\nabla m_k(s)} \leq \delta \norm*{s}}$
        % \Comment{Approximately solve the subproblem}
        \State{%
          $x_{k+1} \gets x_k + s_k$, 
          $r_k \gets \nabla f(x_{k+1}) - \nabla f(x_k) - B_k s_k$
        }
        \State $
          B_{k+1}
          \gets
          \frac{1 - \theta}{1 + \theta} \prn[\Big]{
            B_k
            + \frac{r_k s_k^\top + s_k r_k^\top}{\norm*{s_k}^2}
            - \frac{\inner*{r_k}{s_k}}{\norm*{s_k}^4} s_k s_k^\top
          }
        $
        \label{alg-line:inner_loop_end}
      \EndFor
      \State{%
        $\bar x_K \gets \frac{1}{K (K+1)} \prn[\big]{ \sum_{i=0}^{K-1} (2i+1) x_i + K x_K }$
      }
      \label{alg-line:compute_xbar}
      \If{$\norm*{\nabla f(\bar x_K)}$ is sufficiently small}
        \State \Return $\bar x_K$
      \EndIf
      \State{%
        $x_0 \gets x_K$,
        $B_0 \gets B_K$
      }
      \label{alg-line:prepare_next_outer}
      \Comment{Initialization for the next outer iteration}
    \EndFor
  \end{algorithmic}
\end{algorithm}

\subsection{Update rule for \texorpdfstring{$x_k$}{xk}}
In the inner loop, the iterate $x_k$ is updated using a quasi-Newton method that incorporates a weighted sum of gradients for acceleration together with quartic regularization. 
At inner iteration $k$, we approximately solve the following subproblem to obtain a step $s_k \coloneqq x_{k+1} - x_k$:
\begin{align}
  \!\!
  \min_{s \in \R^d} \, \set[\bigg]{
    m_k(s)
    \coloneqq
    \inner[\Big]{\nabla f(x_k) + \frac{1}{k+1} \sum_{i=0}^k (2i+1) \nabla f(x_i)}{s}
    + \frac{1}{2} \inner*{B_k s}{s}
    + \frac{\sigma}{4} \norm*{s}^4
  }.\!\!\!\!\!
  \label{eq:subproblem}
\end{align}
More precisely, $s_k$ is chosen to satisfy the approximate first-order optimality condition
\begin{align}
  \norm*{\nabla m_k(s_k)}
  \leq
  \delta \norm*{s_k}.
  \label{eq:subproblem_approx_optimality}
\end{align}
Here, $\sigma, \delta > 0$ are fixed throughout the inner loop.
The matrix $B_k \in \SS^d$ approximates the Hessian $\nabla^2 f(x_k)$; its update rule is given in the next section.

Standard regularization methods (e.g., \citep[Algorithm~3.3.1]{cartis2022evaluation}) often impose the model decrease condition $m_k(s_k) < m_k(\0)$ in addition to a first-order condition such as \cref{eq:subproblem_approx_optimality}.
Our analysis avoids relying on a monotonic decrease in the objective value, and therefore allows us to dispense with the model decrease condition.

Let us briefly explain the intuition behind the subproblem~\cref{eq:subproblem}.
Momentum methods such as Nesterov's accelerated gradient~\citep{nesterov1983method} and Polyak's heavy-ball method~\citep{polyak1964methods} are known to accelerate convergence by using a weighted sum of past gradients (e.g., \citep{li2022restarted,li2023restarted,marumo2024parameter,marumo2025universal}).
From another perspective, it is also known that the proximal point method \citep{rockafellar1976monotone}, which updates $x_{k+1} = x_k - \eta_k \nabla f(x_{k+1})$ using the gradient at a \emph{future} point, can converge faster than gradient descent.
The subproblem~\cref{eq:subproblem} is designed by combining these two ideas.
Our aim is to choose the step $s_k$ so that it becomes nearly aligned with $-\bar g_{k+1}$ defined in \cref{eq:def_gbar}, which combines past and future gradient information.
The following lemma shows that if the residual
\begin{align}
  r_k \coloneqq \nabla f(x_{k+1}) - \nabla f(x_k) - B_k s_k
  \label{eq:def_rk}
\end{align}
is sufficiently small, then $s_k$ generated by the subproblem has nearly the same direction as $-\bar g_{k+1}$.

\begin{lemma}\label{lem:subproblem_relaxed_optimality}
  If $s_k \in \R^d$ satisfies the condition \cref{eq:subproblem_approx_optimality}, then the following holds:
  \begin{align}
    \norm*{(k+2) \bar g_{k+1} + \sigma \norm*{s_k}^2 s_k}
    \leq
    \norm*{r_k} + \delta \norm*{s_k}.
    \label{eq:subproblem_relaxed_optimality}
  \end{align}
\end{lemma}
\begin{proof}
  It follows from the definitions of $\bar g_k$ and $r_k$ in \cref{eq:def_rk,eq:def_gbar} that
  \begin{alignat}{2}
    \nabla m_k(s_k)
    &=
    \nabla f(x_k) + \frac{1}{k+1} \sum_{i=0}^k (2i+1) \nabla f(x_i)
    + B_k s_k
    + \sigma \norm*{s_k}^2 s_k\\
    &=
    \nabla f(x_{k+1}) + \frac{1}{k+1} \sum_{i=0}^k (2i+1) \nabla f(x_i)
    - r_k
    + \sigma \norm*{s_k}^2 s_k
    &\quad&\by{\cref{eq:def_rk}}\\
    &=
    (k+2) \bar g_{k+1} - r_k
    + \sigma \norm*{s_k}^2 s_k.
    &\quad&\by{\cref{eq:def_gbar}}
  \end{alignat}
  Combining this with \cref{eq:subproblem_approx_optimality} and applying the triangle inequality completes the proof.
\end{proof}

While quadratic or cubic regularization is more commonly employed for second-order models \citep{ueda2010convergence,mishchenko2023regularized,nesterov2006cubic,doikov2024super}, we adopt quartic regularization because it helps ensure sufficient decrease of the objective value when combined with our update rule for $B_k$.
A detailed explanation is given in the complexity analysis, specifically in \cref{rem:why_quartic}.

\subsection{Update rule for \texorpdfstring{$B_k$}{Bk}}
As shown in the previous section, reducing $\norm*{r_k}$ is desirable for fast convergence.
To achieve this, we update the Hessian approximation $B_k$ using the following scaled variant of the classical Powell-symmetric-Broyden (PSB) update \citep{powell1970new}:
\begin{align}
  B_{k+1}
  =
  \frac{1 - \theta}{1 + \theta} \prn[\bigg]{
    B_k
    + \frac{r_k s_k^\top + s_k r_k^\top}{\norm*{s_k}^2}
    - \frac{\inner*{r_k}{s_k}}{\norm*{s_k}^4} s_k s_k^\top
  },
  \label{eq:update_Bk}
\end{align}
where $\theta \in (0, 1)$ is fixed throughout the inner loop.
As $\theta \to 0$, this formula reduces to the standard PSB update.

To justify the scaled PSB update \cref{eq:update_Bk} and to show that it yields an upper bound on $\norm{r_k}$, let us introduce the following notation:
\begin{align}
  \gamma_k
  &\coloneqq
  \normFro*{B_k - \nabla^2 f(x_k)}^2
  - \normFro*{\nabla^2 f(x_k)}^2
  + \frac{\theta}{1 - \theta} \normFro*{B_k}^2,
  \label{eq:def_gammak}\\
  \quad
  G_k
  &\coloneqq
  \int_0^1 \nabla^2 f(x_k + \tau s_k) \,\dd \tau.
  \label{eq:def_Gk}
\end{align}
The first term in $\gamma_k$ measures the Frobenius distance between $B_k$ and $\nabla^2 f(x_k)$, whereas the second and third terms serve as normalization and regularization, respectively.  
Since $B_k$ is intended to approximate $\nabla^2 f(x_k)$, it is naturally desirable for $\gamma_k$ to be small.  
The matrix $G_k$ represents the averaged Hessian along the segment from $x_k$ to $x_{k+1}$.

The following lemma shows that introducing scaling into quasi-Newton updates, including but not limited to PSB, helps make $\gamma_{k+1}$ small.
The scaling factor $\frac{1 - \theta}{1 + \theta}$ is chosen so that the resulting upper bound on $\gamma_{k+1} - \gamma_k$ takes a simple form.

\begin{lemma}\label{lem:gammak_diff_upperbound_general}
  Suppose that \cref{asm:lipschitz_grad_hess} holds.
  For $\theta \in (0, 1)$ and $B_k, X \in \SS^d$, let $B_{k+1}$ be defined by
  \begin{align}
    B_{k+1}
    &=
    \frac{1 - \theta}{1 + \theta} \prn*{B_k + X}.
    \label{eq:update_Bk_general}
  \end{align}
  Then, the following holds:
  \begin{align}
    \gamma_{k+1} - \gamma_k
    &\leq
    2 d L^2 \theta
    + \frac{d M^2}{2 \theta} \norm*{s_k}^2
    + \normFro*{X}^2 - 2 \inner*{X}{G_k - B_k}.
    \label{eq:gammak_diff_upperbound_general}
  \end{align}
\end{lemma}
The proof is deferred to \cref{sec:proof_lem_gammak_diff_upperbound_general}.

The upper bound~\cref{eq:gammak_diff_upperbound_general} depends on an arbitrary matrix $X \in \SS^d$, so a natural next step is to choose $X$ that minimizes this bound.  
If $X$ could range over the entire space $\SS^d$, the optimal choice would be $X = G_k - B_k$; however, $G_k$, defined in \cref{eq:def_Gk}, is not computable within first-order methods.  
To obtain an implementable update rule for $B_k$, we therefore restrict $X$ to the form $X = v s_k^\top + s_k v^\top$ for some $v \in \R^d$.
As shown in the lemma below, this restriction leads to the scaled PSB update~\cref{eq:update_Bk} and, at the same time, yields an upper bound on $\norm{r_k}$.

\begin{lemma}\label{lem:optimal_X}
  Given $s_k \neq \0$ and $B_k \in \SS^d$, consider the following problem:
  \begin{align}
    \min_{X \in \SS^d} \ 
    \normFro*{X}^2
    - 2 \inner*{X}{G_k - B_k}
    \qquad
    \text{subject to}\ \
    X \in \Set[\big]{
      v s_k^\top + s_k v^\top
    }{
      v \in \R^d
    },
    \label{eq:problem_for_X}
  \end{align}
  where $G_k$ is defined in \cref{eq:def_Gk}.
  The optimal solution is given by
  \begin{align}
    X
    =
    \frac{r_k s_k^\top + s_k r_k^\top}{\norm*{s_k}^2}
    - \frac{\inner*{r_k}{s_k}}{\norm*{s_k}^4} s_k s_k^\top,
    \label{eq:optimal_X}
  \end{align}
  and the optimal value is
  $
    \frac{\inner*{r_k}{s_k}^2}{\norm*{s_k}^4}
    - 2 \frac{\norm*{r_k}^2}{\norm*{s_k}^2}
  $.
\end{lemma}
The proof is deferred to \cref{sec:proof_lem_optimal_X}.

Now, combining \cref{lem:gammak_diff_upperbound_general,lem:optimal_X} shows that the scaled PSB update~\cref{eq:update_Bk} guarantees the following inequality:
\[
  \gamma_{k+1} - \gamma_k
  \leq
  2 d L^2 \theta
  + \frac{d M^2}{2 \theta} \norm*{s_k}^2
  + \frac{\inner*{r_k}{s_k}^2}{\norm*{s_k}^4}
  - 2 \frac{\norm*{r_k}^2}{\norm*{s_k}^2}.
\]
Applying the Cauchy--Schwarz inequality, $\inner*{r_k}{s_k}^2 \leq \norm*{r_k}^2 \norm*{s_k}^2$, and then rearranging terms, we obtain the desired upper bound on $\norm{r_k}$:
\begin{align}
  \frac{\norm*{r_k}^2}{\norm*{s_k}^2}
  \leq
  2 d L^2 \theta
  + \frac{d M^2}{2 \theta} \norm*{s_k}^2
  + \gamma_k - \gamma_{k+1}.
  \label{eq:norm_wk_upperbound}
\end{align}
This bound plays a crucial role in the complexity analysis.

If $s_k = \0$, then $r_k = \0$ by the definition in \cref{eq:def_rk}, and we set $\frac{\norm*{r_k}^2}{\norm*{s_k}^2} = 0$ for convenience.
Under this convention, the bound~\cref{eq:norm_wk_upperbound} holds for all $s_k \in \R^d$.

\subsection{Solving the subproblem}
\label{sec:solving_subproblem}
The subproblem~\cref{eq:subproblem} can be written as
\begin{align}
  \min_{s \in \R^d} \
  \set*{
    m(s)
    \coloneqq
    \inner*{g}{s} + \frac{1}{2} \inner*{Bs}{s} + \frac{\sigma}{4} \norm*{s}^4
  },
  \label{eq:subproblem_general}
\end{align}
where $g \in \R^d$, $B \in \SS^d$, and $\sigma > 0$.
Optimization problems of this form have been extensively studied in the literature; see \citep[Chapters~8--11]{cartis2022evaluation} for a comprehensive review.
Any algorithm can be used as long as it returns $s_k \in \R^d$ satisfying \cref{eq:subproblem_approx_optimality}.
We describe below one approach based on a reduction to a one-dimensional root-finding problem.
The following lemma is a standard result tailored to our setting.

\begin{lemma}\label{lem:subproblem_reduction}
  Suppose that a diagonal matrix $\Lambda = \diag \prn*{\lambda_1, \ldots, \lambda_d}$ and an orthogonal matrix $Q \in \R^{d \times d}$ satisfy  
  \begin{align}
    \normop[\big]{B - Q \Lambda Q^\top}
    \leq
    \frac{\delta}{2},
    \label{eq:approx_eigendecomposition_B}
  \end{align}
  where $\lambda_1 \leq \dots \leq \lambda_d$.
  For $\mu > - \lambda_1$, define
  \begin{align}
    \hat s(\mu) \coloneqq - Q (\Lambda + \mu I)^{-1} Q^\top g,\qquad
    \phi(\mu)
    \coloneqq
    \sigma \norm*{(\Lambda + \mu I)^{-1} Q^\top g}^2
    - \mu.
    \label{eq:def_s_and_phi}
  \end{align}
  Then the following holds:
  \begin{align}
    \norm*{\nabla m(\hat s(\mu))} \leq \prn*{\abs*{\phi(\mu)} + \frac{\delta}{2}} \norm*{\hat s(\mu)}.
  \end{align}
\end{lemma}
\begin{proof}
  Let $\hat B \coloneqq Q \Lambda Q^\top$.
  For notational simplicity, we write $\hat s$ and $\phi$ instead of $\hat s(\mu)$ and $\phi(\mu)$, respectively.
  Since $Q$ is orthogonal, \cref{eq:def_s_and_phi} can be rewritten as
  \begin{align}
    \hat s
    =
    - \prn{\hat B + \mu I}^{-1} g,
    \qquad
    \phi
    =
    \sigma \norm{\hat s}^2 - \mu.
  \end{align}
  Using these equations, we have
  \begin{align}
    g + \hat B \hat s + \sigma \norm*{\hat s}^2 \hat s
    =
    - \prn{\hat B + \mu I} \hat s + \hat B \hat s + \sigma \norm*{\hat s}^2 \hat s
    =
    \prn[\big]{\sigma \norm*{\hat s}^2 - \mu} \hat s
    =
    \phi \hat s.
    \label{eq:subproblem_approx_grad_reduction}
  \end{align}
  Using the triangle inequality, we obtain the desired result:
  \begin{align}
    \norm*{\nabla m(\hat s)}
    =
    \norm[\big]{g + B \hat s + \sigma \norm*{\hat s}^2 \hat s}
    &\leq
    \norm[\big]{g + \hat B \hat s + \sigma \norm*{\hat s}^2 \hat s}
    + \norm[\big]{B - \hat B} \norm*{\hat s}\\
    &\leq
    \norm*{\phi \hat s} + \frac{\delta}{2} \norm*{\hat s},
  \end{align}
  where the last inequality follows from \cref{eq:approx_eigendecomposition_B,eq:subproblem_approx_grad_reduction}.
\end{proof}

The lemma shows that, to satisfy the condition~\cref{eq:subproblem_approx_optimality}, it suffices to find $\mu > - \lambda_1$ such that $\abs*{\phi(\mu)} \leq \delta / 2$.
If $\phi$ has a root, then we can find it using standard root-finding methods such as the bisection method or Newton's method.
Otherwise, a global minimizer of the approximate model obtained by replacing $B$ with $Q \Lambda Q^\top$ can be expressed in closed form; see \citep[Corollary~8.3.1]{cartis2022evaluation} for details.
This minimizer $s \in \R^d$ also satisfies $\norm*{\nabla m(s)} \leq \frac{\delta}{2} \norm{s}$, yielding \cref{eq:subproblem_approx_optimality}.

Let us briefly discuss the computational cost of this approach.
As explained in \cref{sec:approx_eigendecomposition}, the approximate eigendecomposition in \cref{lem:subproblem_reduction} can be computed in $\tilde{\O}(d^2)$ time.
Given $\Lambda$ and $Q$, we first compute $Q^\top g$ in $\O(d^2)$ time.
For each $\mu > - \lambda_1$, evaluating $\phi(\mu)$ takes $\O(d)$ time, since $\Lambda + \mu I$ is diagonal.
The bisection method thus finds a root of $\phi$ to the desired accuracy in $\tilde \O(d)$ time.
We can then compute $\hat s(\mu)$, defined in \cref{eq:def_s_and_phi}, in $\O(d^2)$ time.
If $\phi$ has no root, the global minimizer of the approximate model can also be computed in $\O(d^2)$ time.
Overall, the total cost of solving the subproblem is $\tilde{\O}(d^2)$.

\subsection{Comparison with \texorpdfstring{\citep{jiang2025improved}}{Jiang et al. (2025)}}
\label{sec:comparison_with_jiang}

Although the algorithm of \citep{jiang2025improved} and ours are both based on quasi-Newton methods, there are several important differences between them.
\Cref{table:comparison_with_jiang} summarizes these distinctions.

First, although both methods use an acceleration technique for updating $x_k$, the specific mechanisms are different. 
The method of \citep{jiang2025improved} employs optimistic online learning, whereas ours uses a weighted sum of gradients, similar in spirit to momentum methods. 
While these mechanisms are conceptually related, the resulting algorithms differ in the number of gradient evaluations per iteration: two for \citep{jiang2025improved} versus one for ours.

Second, the subproblem in \citep{jiang2025improved} combines quadratic regularization with an $\ell_2$-norm constraint, and thus it needs two parameters to be specified (the regularization parameter and the constraint radius).
In contrast, our subproblem uses only quartic regularization without any constraints, requiring just one parameter $\sigma$.
This simplification is beneficial for making the algorithm parameter-free.

Third, the update rule for $B_k$ differs significantly.
The method of \citep{jiang2025improved} updates $B_k$ via online learning, and the analysis relies on dynamic regret bounds~\citep[Chapter~14]{orabona2019modern}.
Since standard dynamic regret analyses typically assume that the feasible region is bounded, the previous method employs a separation oracle to keep $\normop*{B_k} \leq 2 L$, which requires prior knowledge of $L$.
In contrast, our update rule for $B_k$ is based on the PSB formula instead of online learning and introduces a simple scaling step.
This scaling step controls the growth of $B_k$ without imposing an explicit norm constraint or requiring knowledge of $L$.

\begin{table}[t]
  \centering
  \caption{Comparison between \citep{jiang2025improved} and \cref{alg:proposed}.}
  \label{table:comparison_with_jiang}
  \small
  \def\arraystretch{1.1}
  \begin{tabular}{@{}l|ll@{}}\toprule
    & \citep{jiang2025improved} & \cref{alg:proposed}\\\midrule
    Acceleration mechanism &
    Optimistic online learning &
    Momentum\\
    Subproblem &
    Quadratic reg.~and $\ell_2$-norm constraint &
    Quartic reg.\\
    $B_k$ update &
    Online learning with separation oracle &
    Scaled PSB\\
    \# Grads./iter. &
    2 &
    1 \\
    Required knowledge &
    $L$, $M$, and $f^*$ &
    None\\
    \bottomrule
  \end{tabular}
\end{table}

\section{Complexity analysis}
\label{sec:complexity_analysis}
We bound the oracle complexity (i.e., the number of function and gradient evaluations required) to reach an $\epsilon$-stationary point for \cref{alg:proposed}.
Since the algorithm evaluates the gradient once per inner iteration, it suffices to bound the total number of inner iterations.

Before proceeding to the main analysis, we present the following lower bound on $\gamma_k$:
\begin{align}
  \gamma_k
  &=
  \frac{1}{1 - \theta} \normFro*{B_k}^2
  - 2 \inner*{B_k}{\nabla^2 f(x_k)}
  \geq
  - (1 - \theta) \normFro*{\nabla^2 f(x_k)}^2
  \geq
  - (1 - \theta) d L^2.
  \label{eq:gamma_equality_lowerbound}
\end{align}
Here, the equality is obtained by expanding the definition of $\gamma_k$ in \cref{eq:def_gammak}, the first inequality follows from Young's inequality, and the second from $\normFro*{\cdot} \leq \sqrt{d} \normop*{\cdot}$ together with $\normop{\nabla^2 f(x_k)} \leq L$ implied by \cref{asm:lipschitz_grad}.
This bound is used in the proofs of \cref{lem:gammak_diff_upperbound_general,lem:min_gradnorm_bound}.

\subsection{Analysis of a single outer iteration}
This section focuses on a single outer iteration, which consists of $K$ inner iterations.
Let us define the potential $\Phi_k$ by
\begin{align}
  \label{eq:def_Phi}
  \Phi_k
  \coloneqq{}&
  f(x_k) + \frac{\gamma_k}{2 \sigma}\\
  ={}&
  f(x_k)
  + \frac{1}{2 \sigma}
  \prn*{
    \normFro*{B_k - \nabla^2 f(x_k)}^2
    - \normFro*{\nabla^2 f(x_k)}^2
    + \frac{\theta}{1 - \theta} \normFro*{B_k}^2
  },
\end{align}
where $\gamma_k$ is defined in \cref{eq:def_gammak}.
The goal of this section is to prove the following lemma, which bounds the gradient norm at the averaged solution $\bar x_K$ in terms of the potential decrease $\Phi_0 - \Phi_K$.

\begin{lemma}\label{lem:gradnorm_upperbound_by_potential_decrease}
  Suppose that \cref{asm:lipschitz_grad_hess} holds.
  Let $\kappa > d^{1/5}$, $\theta = d / \kappa^5$, and $K = \floor*{\kappa}$.
  For the sequence $(x_k)_{k=0}^K$ generated by the inner loop of \cref{alg:proposed}, the following holds:
  \[
    \frac{M \kappa^4}{\sigma}
    \norm*{\nabla f(\bar x_K)}
    \leq
    \Phi_0 - \Phi_K
    + \frac{d^2 L^2}{\sigma \kappa^4}
    + \frac{\delta^2 \kappa}{2 \sigma}
    + \frac{5 M^4 \kappa^{12}}{2 \sigma^3},
  \]
  where $\bar x_K$ is defined by \cref{eq:def_xbar}.
\end{lemma}

% The lemma above can also be viewed as providing a lower bound on the potential decrease.
% Unlike vanilla gradient descent, it is difficult to guarantee a lower bound on the objective decrease itself for the proposed algorithm.
% However, the potential $\Phi_k$, which is defined using both the objective value and the Hessian approximation $B_k$, does admit such a guarantee.
% Roughly speaking, after $K$ iterations, either the objective value or $\gamma_k$, which primarily measures the mismatch between $B_k$ and the true Hessian, must decrease sufficiently.

To prove \cref{lem:gradnorm_upperbound_by_potential_decrease}, we first establish two auxiliary lemmas.
The first lemma relates the gradient at the averaged solution $\bar x_k$ to the averaged gradient $\bar g_k$.
To simplify the notation, let
\begin{align}
  S_k
  &\coloneqq
  \sum_{i=0}^{k-1} \norm*{s_i}^2.
  \label{eq:def_Sk}
\end{align}
\begin{lemma}\label{lem:grad_norm_at_average}
  Suppose that \cref{asm:lipschitz_hess} holds.
  Let $(x_k)_{k \in \N}$ be the sequence generated by the inner loop of \cref{alg:proposed}.
  Then, the following holds for all $k \geq 1$:
  \begin{align}
    \norm*{\nabla f(\bar x_k)}
    \leq
    \norm*{\bar g_k}
    + \frac{M}{8} k S_k,
    \label{eq:grad_norm_at_average}
  \end{align}
  where $\bar x_k$ and $\bar g_k$ are defined in \cref{eq:def_xbar,eq:def_gbar}, respectively.
\end{lemma}
The proof is deferred to \cref{sec:proof_lem_grad_norm_at_average} as it is almost a direct application of a Hessian-free inequality~\citep[Lemma~3.1]{marumo2024parameter}.
The Hessian-free inequality is equivalent to \cref{asm:lipschitz_hess} \citep{bot2025equivalence}.

The second lemma is more involved, providing an upper bound on $f(x_k)$.
\begin{lemma}\label{lem:function_decrease}
  Suppose that \cref{asm:lipschitz_hess} holds.
  Let $(x_k)_{k \in \N}$ be the sequence generated by the inner loop of \cref{alg:proposed}.
  Then, the following holds for all $k \geq 1$:
  \begin{align}
    f(x_k) - f(x_0)
    &\leq
    \frac{\delta^2 k}{2 \sigma}
    - \frac{3}{8} \frac{\norm*{(k+1) \bar g_k}^{4/3}}{\sigma^{1/3}}
    - \frac{\sigma}{8 k^2} S_k^2
    + \frac{M}{12} S_k^{3/2}
    + \frac{1}{2 \sigma} \sum_{i=0}^{k-1} \frac{\norm*{r_i}^2}{\norm*{s_i}^2},
    \label{eq:function_decrease}
  \end{align}
  where we set $\frac{\norm*{r_i}^2}{\norm*{s_i}^2} = 0$ if $s_i = \0$.
\end{lemma}
\begin{proof}
  A Hessian-free inequality in \citep[Lemma~3.2]{marumo2024parameter} gives
  \begin{align}
    f(x_{i+1}) - f(x_i)
    \leq
    \frac{1}{2} \inner*{\nabla f(x_i) + \nabla f(x_{i+1})}{s_i}
    + \frac{M}{12} \norm*{s_i}^3.
    \label{eq:hessian-free-ineq-function}
  \end{align}
  Summing this inequality over $i = 0,\dots,k-1$ yields
  \begin{align}
    f(x_k) - f(x_0)
    &\leq
    \frac{1}{2} \sum_{i=0}^{k-1} \inner*{\nabla f(x_i) + \nabla f(x_{i+1})}{s_i}
    + \frac{M}{12} \sum_{i=0}^{k-1} \norm*{s_i}^3.
    \label{eq:func_dec_proof1}
  \end{align}
  \Cref{eq:relation_gbar_g} allows us to rewrite the first sum on the right-hand side as
  \begin{align}
    \sum_{i=0}^{k-1} \inner*{\nabla f(x_i) + \nabla f(x_{i+1})}{s_i}
    &=
    \sum_{i=0}^{k-1} (i+2) \inner*{\bar g_{i+1}}{s_i}
    - \sum_{i=1}^{k-1} i \inner*{\bar g_i}{s_i}.
    \label{eq:func_dec_proof2}
  \end{align}
  We now bound the two terms in \cref{eq:func_dec_proof2} separately.

  For the first term, we square both sides of \cref{eq:subproblem_relaxed_optimality} to obtain
  \begin{align}
    (i+2)^2 \norm*{\bar g_{i+1}}^2
    + 2 (i+2) \sigma \norm*{s_i}^2 \inner*{\bar g_{i+1}}{s_i}
    + \sigma^2 \norm*{s_i}^6
    &\leq
    \prn*{
      \norm*{r_i} + \delta \norm*{s_i}
    }^2\\
    &\leq
    2 \norm*{r_i}^2
    + 2 \delta^2 \norm*{s_i}^2,
  \end{align}
  where $(a + b)^2 \leq 2 (a^2 + b^2)$ is used.
  Dividing both sides by $2 \sigma \norm*{s_i}^2$ and rearranging terms yields
  \begin{align}
    \label{eq:upperbound_inner_gbar_s_1}
    \qquad
    (i+2) \inner*{\bar g_{i+1}}{s_i}
    &\leq
    \frac{\delta^2}{\sigma}
    + \frac{\norm*{r_i}^2}{\sigma \norm*{s_i}^2}
    - \frac{\sigma}{2} \norm*{s_i}^4
    - \frac{\norm*{(i+2) \bar g_{i+1}}^2}{2 \sigma \norm*{s_i}^2}\\
    &=
    \frac{\delta^2}{\sigma}
    + \frac{\norm*{r_i}^2}{\sigma \norm*{s_i}^2}
    - \frac{\sigma}{4} \norm*{s_i}^4
    - \frac{3}{4} \prn*{
      \frac{\sigma}{3} \norm*{s_i}^4
      + \frac{2}{3} \frac{\norm*{(i+2) \bar g_{i+1}}^2}{\sigma \norm*{s_i}^2}
    }\\
    &\leq
    \frac{\delta^2}{\sigma}
    + \frac{\norm*{r_i}^2}{\sigma \norm*{s_i}^2}
    - \frac{\sigma}{4} \norm*{s_i}^4
    - \frac{3}{4} 
    \prn*{
      \sigma \norm*{s_i}^4
    }^{1/3}
    \prn[\bigg]{
      \frac{\norm*{(i+2) \bar g_{i+1}}^2}{\sigma \norm*{s_i}^2}
    }^{2/3}\\
    &=
    \frac{\delta^2}{\sigma}
    + \frac{\norm*{r_i}^2}{\sigma \norm*{s_i}^2}
    - \frac{\sigma}{4} \norm*{s_i}^4
    - \frac{3}{4} \frac{\norm*{(i+2) \bar g_{i+1}}^{4/3}}{\sigma^{1/3}},
  \end{align}
  where the last inequality follows from the weighted AM--GM inequality $\frac{1}{3} a + \frac{2}{3} b \geq a^{1/3} b^{2/3}$.
  Note that \cref{eq:upperbound_inner_gbar_s_1} is valid even when $s_i = \0$, because \cref{eq:subproblem_relaxed_optimality} then implies $\bar g_{i+1} = \0$.

  For the second term on the right-hand side of \cref{eq:func_dec_proof2}, we again use weighted AM--GM $a^{1/4} b^{3/4} \leq \frac{1}{4}a + \frac{3}{4}b$ to obtain
  \begin{align}
    - i \inner*{\bar g_i}{s_i}
    \leq
    i \norm*{\bar g_i} \norm*{s_i}
    &=
    \prn[\bigg]{
      \frac{\sigma i^4 \norm*{s_i}^4}{(i+1)^4}
    }^{1/4}
    \prn[\bigg]{
      \frac{\norm*{(i+1) \bar g_i}^{4/3}}{\sigma^{1/3}}
    }^{3/4}\\
    &\leq
    \frac{1}{4} \frac{\sigma i^4 \norm*{s_i}^4}{(i+1)^4}
    + \frac{3}{4} \frac{\norm*{(i+1) \bar g_i}^{4/3}}{\sigma^{1/3}}.
  \end{align}
  Substituting this bound and \cref{eq:upperbound_inner_gbar_s_1} into \cref{eq:func_dec_proof2} and noting that the intermediate terms involving $\norm*{\bar g_i}$ cancel out yields
  \begin{align}
    &\mathInd
    {\sum_{i=0}^{k-1}} \inner*{\nabla f(x_i) + \nabla f(x_{i+1})}{s_i}\\
    &\leq
    \sum_{i=0}^{k-1} \prn[\bigg]{
      \frac{\delta^2}{\sigma}
      + \frac{\norm*{r_i}^2}{\sigma \norm*{s_i}^2}
      - \frac{\sigma}{4} \norm*{s_i}^4
      - \frac{3}{4} \frac{\norm*{(i+2) \bar g_{i+1}}^{4/3}}{\sigma^{1/3}}
    }\\
    &\qquad
    + \sum_{i=1}^{k-1} \prn*{
      \frac{1}{4} \frac{\sigma i^4 \norm*{s_i}^4}{(i+1)^4}
      + \frac{3}{4} \frac{\norm*{(i+1) \bar g_i}^{4/3}}{\sigma^{1/3}}
    }\\
    &=
    \frac{\delta^2 k}{\sigma}
    + \sum_{i=0}^{k-1} \frac{\norm*{r_i}^2}{\sigma \norm*{s_i}^2}
    + \frac{\sigma}{4} \sum_{i=0}^{k-1} \prn*{ \frac{i^4}{(i+1)^4} - 1 } \norm*{s_i}^4
    - \frac{3}{4} \frac{\norm*{(k+1) \bar g_k}^{4/3}}{\sigma^{1/3}}\\
    &\leq
    \frac{\delta^2 k}{\sigma}
    + \sum_{i=0}^{k-1} \frac{\norm*{r_i}^2}{\sigma \norm*{s_i}^2}
    - \frac{\sigma}{4k} \sum_{i=0}^{k-1} \norm*{s_i}^4
    - \frac{3}{4} \frac{\norm*{(k+1) \bar g_k}^{4/3}}{\sigma^{1/3}},
    \label{eq:func_dec_proof3}
  \end{align}
  where the last inequality follows from $\frac{i^4}{(i+1)^4} \leq \frac{i}{i+1} \leq \frac{k-1}{k}$ for $1 \leq i \leq k-1$.
  Plugging this bound into \cref{eq:func_dec_proof1} yields
  \begin{align}
    \label{eq:func_dec_proof4}
    &\mathInd
    f(x_k) - f(x_0)\\
    &\leq
    \frac{\delta^2 k}{2 \sigma}
    + \frac{1}{2 \sigma} \sum_{i=0}^{k-1} \frac{\norm*{r_i}^2}{\norm*{s_i}^2}
    - \frac{3}{8} \frac{\norm*{(k+1) \bar g_k}^{4/3}}{\sigma^{1/3}}
    - \frac{\sigma}{8 k} \sum_{i=0}^{k-1} \norm*{s_i}^4
    + \frac{M}{12} \sum_{i=0}^{k-1} \norm*{s_i}^3.
  \end{align}

  Finally, we bound the last two terms.
  The Cauchy--Schwarz inequality gives
  \[
    S_k^2
    = 
    \prn*{\sum_{i=0}^{k-1} \norm*{s_i}^2}^2
    \leq
    \prn*{\sum_{i=0}^{k-1} 1^2} \prn*{\sum_{i=0}^{k-1} \norm*{s_i}^4}
    =
    k \sum_{i=0}^{k-1} \norm*{s_i}^4.
  \]
  The monotonicity of $\ell_p$-norms gives $\norm*{a}_3 \leq \norm*{a}_2$ for any $a \in \R^k$, and hence $\sum_{i=0}^{k-1} \norm*{s_i}^3 \allowbreak \leq \prn[\big]{\sum_{i=0}^{k-1} \norm*{s_i}^2}^{3/2} = S_k^{3/2}$.
  Plugging these bounds into \cref{eq:func_dec_proof4} completes the proof.
\end{proof}

We apply \cref{eq:norm_wk_upperbound} to bound the term $\frac{\norm*{r_i}^2}{\norm*{s_i}^2}$ appearing in \cref{eq:function_decrease}.
Combining the resulting bound with \cref{lem:grad_norm_at_average}, we prove \cref{lem:gradnorm_upperbound_by_potential_decrease}.

\begin{remark}\label{rem:why_quartic}
  The appearance of the same term $\frac{\norm*{r_i}^2}{\norm*{s_i}^2}$ in both \cref{eq:norm_wk_upperbound,eq:function_decrease} is no accident, but a design choice.
  Indeed, when the subproblem~\cref{eq:subproblem} employs $p$th-order regularization, we can still prove an analogue of \cref{eq:function_decrease}, yet the resulting bound involves $\frac{\norm*{r_i}^2}{\norm*{s_i}^{p-2}}$.
  We therefore set $p=4$ so that this quantity matches $\frac{\norm*{r_i}^2}{\norm*{s_i}^2}$ bounded in \cref{eq:norm_wk_upperbound}.
  This is the principal reason we adopt quartic regularization for the subproblem.
\end{remark}

We now present the proof of \cref{lem:gradnorm_upperbound_by_potential_decrease}.

\begin{proof}[Proof of \cref{lem:gradnorm_upperbound_by_potential_decrease}]
  Summing \cref{eq:norm_wk_upperbound} over $k$ yields
  \begin{align}
    \sum_{i=0}^{K-1} \frac{\norm*{r_i}^2}{\norm*{s_i}^2}
    &\leq
    \sum_{i=0}^{K-1} \prn*{
      2 d L^2 \theta
      + \frac{d M^2}{2 \theta} \norm*{s_i}^2
      + \gamma_i - \gamma_{i+1}
    }\\
    &=
    2 d L^2 \theta K
    + \frac{d M^2}{2 \theta} S_K
    + \gamma_0 - \gamma_K
    =
    \frac{2 d^2 L^2 K}{\kappa^5}
    + \frac{M^2 \kappa^5}{2} S_K
    + \gamma_0 - \gamma_K,
  \end{align}
  where the last equality uses $\theta = d / \kappa^5$.
  Substituting this bound into \cref{lem:function_decrease} and recalling the definition of $\Phi_k$ in \cref{eq:def_Phi}, we obtain
  \begin{align}
    0
    &\leq
    f(x_0) - f(x_K)
    + \frac{\delta^2 K}{2 \sigma}
    - \frac{3}{8} \frac{\norm*{(K+1) \bar g_K}^{4/3}}{\sigma^{1/3}}
    - \frac{\sigma}{8 K^2} S_K^2
    + \frac{M}{12} S_K^{3/2}\\
    &\qquad
    + \frac{1}{2 \sigma} \prn*{
      \frac{2 d^2 L^2 K}{\kappa^5}
      + \frac{M^2 \kappa^5}{2} S_K
      + \gamma_0 - \gamma_K
    }\\
    &=
    \Phi_0 - \Phi_K
    + \prn*{
      \frac{d^2 L^2}{\kappa^5}
      + \frac{\delta^2}{2}
    } \frac{K}{\sigma}
    - \frac{3}{8} \frac{\norm*{(K+1) \bar g_K}^{4/3}}{\sigma^{1/3}}\\
    &\qquad
    + \frac{M^2 \kappa^5}{4 \sigma} S_K
    + \frac{M}{12} S_K^{3/2}
    - \frac{\sigma}{8 K^2} S_K^2.
  \end{align}
  Using $K \leq \kappa < K + 1$, we replace $K$ and $(K+1)$ in the coefficients by $\kappa$, and we multiply the term $\frac{M}{12} S_K^{3/2}$ by $\kappa^{3/2} \geq 1$, yielding
  \[
  % \begin{align}
    0
    \leq
    \Phi_0 - \Phi_K
    + \prn*{
      \frac{d^2 L^2}{\sigma \kappa^4}
      + \frac{\delta^2 \kappa}{2 \sigma}
    }
    - \frac{3}{8} \frac{\norm*{\kappa \bar g_K}^{4/3}}{\sigma^{1/3}}
    + \frac{M^2 \kappa^5}{4 \sigma} S_K
    + \frac{M}{12} (\kappa S_K)^{3/2}
    - \frac{\sigma}{8 \kappa^2} S_K^2.
    % \label{eq:grad_norm_upperbound_proof0}
  % \end{align}
  \]
  \cref{lem:grad_norm_at_average} gives
  \[
    \norm*{\nabla f(\bar x_K)}
    \leq
    \norm*{\bar g_K}
    + \frac{M}{8} K S_K
    \leq
    \norm*{\bar g_K}
    + \frac{M}{8} \kappa S_K.
  \]
  Multiplying both sides by $\frac{M \kappa^4}{\sigma}$ and adding it to the previous bound gives
  \begin{align}
    \label{eq:grad_norm_upperbound_proof1}
    \quad\ 
    \frac{M \kappa^4}{\sigma} \norm*{\nabla f(\bar x_K)}
    &\leq
    \Phi_0 - \Phi_K
    + \prn*{
      \frac{d^2 L^2}{\sigma \kappa^4}
      + \frac{\delta^2 \kappa}{2 \sigma}
    }
    + \overbrace{
      \prn[\bigg]{
        \frac{M \kappa^4}{\sigma} \norm*{\bar g_K}
        - \frac{3}{8} \frac{\norm*{\kappa \bar g_K}^{4/3}}{\sigma^{1/3}}
      }
    }^{\text{(A)}}\!\!\!\\
    &\qquad
    + \underbrace{
      \prn*{
        \frac{3}{8} \frac{M^2 \kappa^5}{\sigma} S_K
        + \frac{M}{12} (\kappa S_K)^{3/2}
        - \frac{\sigma}{8 \kappa^2} S_K^2
      }
    }_{\text{(B)}}.
  \end{align}

  We next bound the terms (A) and (B) separately.
  To bound (A), we apply the weighted AM--GM inequality $\frac{3}{4} a + \frac{1}{4} b \geq a^{3/4} b^{1/4}$:
  \[
    \frac{3}{4} \prn[\bigg]{\frac{\norm*{\kappa \bar g_K}^{4/3}}{2 \sigma^{1/3}}}
    + \frac{1}{4} \prn*{8 \frac{M^4 \kappa^{12}}{\sigma^3}}
    \geq
    \prn[\bigg]{\frac{\norm*{\kappa \bar g_K}^{4/3}}{2 \sigma^{1/3}}}^{3/4}
    \prn*{8 \frac{M^4 \kappa^{12}}{\sigma^3}}^{1/4}
    =
    \frac{M \kappa^4}{\sigma} \norm*{\bar g_K}.
  \]
  Rearranging terms yields an upper bound on (A): $
    \text{(A)}
    \leq
    2 \frac{M^4 \kappa^{12}}{\sigma^3}$.
  To bound (B) in \cref{eq:grad_norm_upperbound_proof1}, we rewrite it as
  \[
    \text{(B)}
    =
    \frac{M^4 \kappa^{12}}{\sigma^3}
    \prn*{
      \frac{3}{8} \xi^2 
      + \frac{1}{12} \xi^3
      - \frac{1}{8} \xi^4
    },
    \quad\text{where}\quad
    \xi
    \coloneqq
    \frac{\sigma}{M} \sqrt{\frac{S_K}{\kappa^7}}.
  \]
  Let $\phi(\xi) \coloneqq \frac{3}{8} \xi^2 + \frac{1}{12} \xi^3 - \frac{1}{8} \xi^4$.
  Since $\phi'(\xi) = \frac{3}{4} \xi + \frac{1}{4} \xi^2 - \frac{1}{2} \xi^3 = - \frac{1}{4} \xi (\xi + 1) (2 \xi - 3)$, the unique maximizer of $\phi(\xi)$ over $\xi \geq 0$ is $\xi = 3/2$, and hence $\phi(\xi) \leq \phi(\frac{3}{2}) = \frac{63}{128} < \frac{1}{2}$.
  Therefore, (B) can be bounded as
  $
    \text{(B)}
    \leq
    \frac{1}{2} \frac{M^4 \kappa^{12}}{\sigma^3}
  $.

  We now have $\text{(A)} + \text{(B)} \leq \frac{5}{2} \frac{M^4 \kappa^{12}}{\sigma^3}$.
  Plugging this bound into \cref{eq:grad_norm_upperbound_proof1} completes the proof.  
\end{proof}

\subsection{General bound on the gradient norm}
The previous section analyzed a single outer iteration, while this section analyzes the entire algorithm to derive an oracle complexity bound.
Note that parameters such as $\sigma$ and $\kappa$ vary across outer iterations.

Let $\kappa_t$, $\sigma_t$, $\delta_t$, $\theta_t$, and $K_t$ denote the values of $\kappa$, $\sigma$, $\delta$, $\theta$, and $K$ in outer iteration $t \in \N$.
Similarly, we write $x_k$, $\bar{x}_k$, $B_k$, $\gamma_k$, and $\Phi_k$ as $x_k^t$, $\bar{x}_k^t$, $B_k^t$, $\gamma_k^t$, and $\Phi_k^t$, respectively.
With this notation, the upper bound in \cref{lem:gradnorm_upperbound_by_potential_decrease} can be rewritten as
\begin{align}
  \frac{M \kappa_t^4}{\sigma_t}
  \norm*{\nabla f(\bar x_{K_t}^t)}
  &\leq
  \Phi^t_0 - \Phi^t_{K_t}
  + \frac{d^2 L^2}{\sigma_t \kappa_t^4}
  + \frac{\delta_t^2 \kappa_t}{2 \sigma_t}
  + \frac{5 M^4 \kappa_t^{12}}{2 \sigma_t^3}.
  \label{eq:grad_norm_upperbound_t}
\end{align}
Using this inequality, we first derive the following general bound on the gradient norm.
Let
\[
  \Delta \coloneqq f(\xinit) - \inf_{x \in \R^d} f(x),
\]
where $\xinit \in \R^d$ is the initial point used in \cref{alg:proposed}.
\begin{lemma}
  \label{lem:min_gradnorm_bound}
  Suppose that \cref{asm:lipschitz_grad_hess} holds and that the sequences $(\sigma_t)_{t \in \N}$ and $(\kappa_t)_{t \in \N}$ are nondecreasing.
  Then, the following holds for all $T \geq 1$:
  \begin{align}
    \min_{0 \leq t < T}
    \norm*{\nabla f(\bar x_{K_t}^t)}
    &\leq
    \prn*{
      \sum_{t=0}^{T-1} \frac{M \kappa_t^4}{\sigma_t}
    }^{-1}
    \prn*{
      \Delta + \frac{dL^2}{2 \sigma_0}
      +
      \sum_{t=0}^{T-1}
      \prn*{
        \frac{d^2 L^2}{\sigma_t \kappa_t^4}
        + \frac{\delta_t^2 \kappa_t}{2 \sigma_t}
        + \frac{5 M^4 \kappa_t^{12}}{2 \sigma_t^3}
      }
    }.
    \label{eq:gradnorm_bound_general}
  \end{align}  
\end{lemma}
\begin{proof}
  Summing \cref{eq:grad_norm_upperbound_t} over $t = 0,\dots,T-1$ gives
  \begin{align}
    \label{eq:overall_gradnorm_proof0}
    &\mathInd
    \prn*{
      \sum_{t=0}^{T-1} \frac{M \kappa_t^4}{\sigma_t}
    }
    \min_{0 \leq t < T} \norm*{\nabla f(\bar x^t_{K_t})}\\
    &\leq
    \sum_{t=0}^{T-1} \frac{M \kappa_t^4}{\sigma_t}
    \norm*{\nabla f(\bar x^t_{K_t})}\\
    &\leq
    \sum_{t=0}^{T-1} \prn*{
      \Phi^t_0 - \Phi^t_{K_t}
    }
    + \sum_{t=0}^{T-1} \prn*{
      \frac{d^2 L^2}{\sigma_t \kappa_t^4}
      + \frac{\delta_t^2 \kappa_t}{2 \sigma_t}
      + \frac{5 M^4 \kappa_t^{12}}{2 \sigma_t^3}
    }.
  \end{align}
  To bound the first sum on the right-hand side, we decompose it as follows:
  \begin{align}
    \sum_{t=0}^{T-1} \prn*{ \Phi^t_0 - \Phi^t_{K_t} }
    &=
    \Phi^0_0
    - \Phi^T_0
    + \sum_{t=1}^T \prn*{ \Phi^{t}_0 - \Phi^{t-1}_{K_{t-1}} }.
    \label{eq:potential_decrease_decomposition}
  \end{align}

  For the first term, since $x^0_0 = \xinit$ and $B^0_0 = O$, we have
  \begin{align}
    \Phi^0_0
    &=
    f(x^0_0)
    + \frac{1}{2 \sigma_0} \prn*{
      \normFro*{B^0_0 - \nabla^2 f(x^0_0)}^2
      - \normFro*{\nabla^2 f(x^0_0)}^2
      + \frac{\theta_0}{1 - \theta_0} \normFro*{B^0_0}^2
    }
    =
    f(\xinit).
    \label{eq:potential_decrease_decomposition_term1}
  \end{align}
  For the second term, we have
  \begin{align}
    - \Phi^T_0
    &=
    - f(x^T_0)
    - \frac{\gamma_0^T}{2 \sigma_T}
    \leq
    - f(x^T_0)
    + \frac{(1 - \theta_T) dL^2}{2 \sigma_T}
    \leq
    - f(x^T_0)
    + \frac{d L^2}{2 \sigma_T},
    \label{eq:potential_decrease_decomposition_term2}
  \end{align}
  where the first inequality follows from \cref{eq:gamma_equality_lowerbound}.
  For the last term in \cref{eq:potential_decrease_decomposition}, note that $x^{t-1}_{K_{t-1}} = x^t_0$ and $B^{t-1}_{K_{t-1}} = B^t_0$ by Line~\ref{alg-line:prepare_next_outer} of \cref{alg:proposed}.
  Therefore, we have
  \begin{align}
    &\mathInd
    \Phi_0^t - \Phi_{K_{t-1}}^{t-1}\\
    &=
    \prn*{
      f(x^t_0)
      + \frac{1}{2 \sigma_{t}} \prn*{
        \normFro*{B^t_0 - \nabla^2 f(x^t_0)}^2
        - \normFro*{\nabla^2 f(x^t_0)}^2
        + \frac{\theta_t}{1 - \theta_t} \normFro*{B^t_0}^2
      }
    }\\
    &\qquad
    - \prn*{
      f(x^t_0)
      + \frac{1}{2 \sigma_{t-1}} \prn*{
        \normFro*{B^t_0 - \nabla^2 f(x^t_0)}^2
        - \normFro*{\nabla^2 f(x^t_0)}^2
        + \frac{\theta_{t-1}}{1 - \theta_{t-1}} \normFro*{B^t_0}^2
      }
    }\\
    &=
    \frac{1}{2} \prn*{ \frac{1}{\sigma_{t-1}} - \frac{1}{\sigma_t} } \prn*{
      \normFro*{\nabla^2 f(x_0^t)}^2
      - \normFro*{B_{0}^t - \nabla^2 f(x_{0}^t)}^2
    }\\
    &\qquad
    + \frac{1}{2} \prn*{
      \frac{\theta_{t}}{\sigma_{t} (1 - \theta_{t})}
      - \frac{\theta_{t-1}}{\sigma_{t-1} (1 - \theta_{t-1})}
    } \normFro*{B_{0}^t}^2.
  \end{align}
  Since $\sigma_t$ is nondecreasing and $\theta_t = d / \kappa_t^5$ is nonincreasing, the coefficients satisfy
  \[
    \frac{1}{\sigma_{t-1}} - \frac{1}{\sigma_t}
    \geq 0,
    \quad
    \frac{\theta_{t}}{\sigma_{t} (1 - \theta_{t})}
    - \frac{\theta_{t-1}}{\sigma_{t-1} (1 - \theta_{t-1})}
    \leq 0.
  \]
  Hence, we have
  \begin{align}
    \Phi_0^t - \Phi_{K_{t-1}}^{t-1}
    &\leq
    \frac{1}{2} \prn*{ \frac{1}{\sigma_{t-1}} - \frac{1}{\sigma_t} }
    \normFro*{\nabla^2 f(x_{0}^t)}^2
    \leq
    \frac{dL^2}{2} \prn*{ \frac{1}{\sigma_{t-1}} - \frac{1}{\sigma_t} }.
    \label{eq:potential_decrease_decomposition_term3}
  \end{align}
  Plugging \cref{eq:potential_decrease_decomposition_term1,eq:potential_decrease_decomposition_term2,eq:potential_decrease_decomposition_term3} into \cref{eq:potential_decrease_decomposition} yields
  \begin{align}
    \sum_{t=0}^{T-1} \prn*{ \Phi^t_0 - \Phi^t_{K_t} }
    &\leq
    f(\xinit) - f(x^T_0)
    + \frac{d L^2}{2 \sigma_T}
    + \sum_{t=1}^T \frac{dL^2}{2} \prn*{ \frac{1}{\sigma_{t-1}} - \frac{1}{\sigma_t} }\\
    &=
    f(\xinit) - f(x^T_0) + \frac{dL^2}{2 \sigma_0}
    \leq
    \Delta + \frac{dL^2}{2 \sigma_0}.
  \end{align}
  Combining this bound with \cref{eq:overall_gradnorm_proof0} completes the proof.
\end{proof}

\subsection{Oracle complexity}
Let us first analyze the oracle complexity of \cref{alg:proposed} in the case where the number of outer iterations $T$ is given in advance.
To balance the summation terms appearing in \cref{eq:gradnorm_bound_general}, we then set the parameters as
\begin{align}
  \kappa_t
  = c_\kappa T^{1/12},\quad
  \sigma_t
  = c_\sigma T^{2/3},\quad
  \delta_t
  = c_\delta T^{-5/24},
  \label{eq:parameter_schedule_given_T}
\end{align}
where $c_\kappa > d^{1/5}$ and $c_\sigma, c_\delta > 0$ are constants.
The derivation of these exponents is given in \cref{sec:parameter_optimality}.
This parameter setting yields the following oracle complexity bound.
\begin{theorem}
  \label{thm:oracle_complexity_given_T}
  Suppose that \cref{asm:lipschitz_grad_hess} holds.
  Set the parameters $\kappa_t$, $\sigma_t$, and $\delta_t$ as in \cref{eq:parameter_schedule_given_T}.
  Then the following holds for all $T \geq 1$:
  \begin{align}
    \min_{0 \leq t < T}
    \norm*{\nabla f(\bar x_{K_t}^t)}
    \leq
    \frac{C}{T^{2/3}}
    + \frac{dL^2}{2 M c_\kappa^4 T^{4/3}},
    \label{eq:gradnorm_bound_given_T}
  \end{align}
  where
  \begin{align}
    C
    \coloneqq
    \frac{1}{M}
    \prn*{
      \frac{\Delta c_\sigma}{c_\kappa^4}
      + \frac{d^2 L^2}{c_\kappa^8}
      + \frac{c_\delta^2}{2 c_\kappa^3}
      + \frac{5 M^4 c_\kappa^8}{2 c_\sigma^2}
    }.
    \label{eq:def_C}
  \end{align}
  Furthermore, with this parameter setting, \cref{alg:proposed} finds an $\epsilon$-stationary point within
  \begin{align}
    \O \prn*{
      c_\kappa \prn*{\frac{C}{\epsilon}}^{13/8}
      + c_\kappa \prn*{\frac{dL^2}{M c_\kappa^4 \epsilon}}^{13/16}
    }
    \label{eq:oracle_complexity_given_T}
  \end{align}
  gradient evaluations.
\end{theorem}
\begin{proof}
  Plugging \cref{eq:parameter_schedule_given_T} into \cref{eq:gradnorm_bound_general} proves \cref{eq:gradnorm_bound_given_T}.
  The inequality \cref{eq:gradnorm_bound_given_T} implies that
  \begin{align}
    T
    =
    \O \prn*{
      \prn*{\frac{C}{\epsilon}}^{3/2}
      + \prn*{\frac{dL^2}{M c_\kappa^4 \epsilon}}^{3/4}
    }
    \label{eq:num_outer_iterations}
  \end{align}
  suffices to find an $\epsilon$-stationary point.
  Then, the total number of inner iterations is at most
  \[
    \sum_{t=0}^{T-1} \kappa_t
    =
    c_\kappa T^{13/12}
    =
    \O \prn*{
      c_\kappa \prn*{\frac{C}{\epsilon}}^{13/8}
      + c_\kappa \prn*{\frac{dL^2}{M c_\kappa^4 \epsilon}}^{13/16}
    },
  \]
  which completes the proof.
\end{proof}

We now discuss how to choose the constants $c_\kappa$, $c_\sigma$, and $c_\delta$.
Since the problem dimension $d$ is typically known, we set
\begin{align}
  c_\kappa
  =
  \Theta \prn[\big]{ d^{1/4} },\quad
  c_\sigma
  =
  \Theta \prn*{ d },\quad
  c_\delta
  =
  \Theta \prn[\big]{ d^{3/8} }.
  \label{eq:def_cs_only_d}
\end{align}
With these choices, the constant $C$ defined in \cref{eq:def_C} becomes
$C = \Theta \prn[\big]{\frac{\Delta + L^2 + 1}{M} + M^3}$.
Hence, the oracle complexity \cref{eq:oracle_complexity_given_T} becomes
% \begin{align}
\begin{align}
  \O \prn*{
    d^{1/4}
    \prn*{
      \frac{\Delta + L^2 + 1}{M \epsilon} + \frac{M^3}{\epsilon}
    }^{13/8}
    + d^{1/4} \prn*{\frac{L^2}{M \epsilon}}^{13/16}
  }
  =
  \O \prn[\bigg]{
    \frac{d^{1/4}}{\epsilon^{13/8}}
  },
  \label{eq:oracle_complexity_only_d}
\end{align}
where the last $\O$-notation ignores dependencies on parameters other than $d$ and $\epsilon$.

If the scales of $L$, $M$, and $\Delta$ are also known, the dependence can be improved by setting
\begin{align}
  c_\kappa
  =
  \Theta \prn[\bigg]{ \frac{(d L)^{1/4}}{(M \sqrt{\Delta})^{1/6}} },\quad
  c_\sigma
  =
  \Theta \prn*{ \frac{d L M^{2/3}}{\Delta^{2/3}} },\quad
  c_\delta
  =
  \Theta \prn*{ (d L)^{3/8} (M \sqrt{\Delta})^{5/12} }.
  \label{eq:setting_cs_with_d_L_M_Delta}
\end{align}
In this case, we have $C = \Theta \prn[\big]{M^{1/3} \Delta^{2/3}}$, and the oracle complexity becomes
\begin{align}
  \O \prn*{
    \frac{(d L)^{1/4} M^{3/8} \Delta}{\epsilon^{13/8}}
  }
  + \O \prn*{\frac{d^{1/4}}{\epsilon^{13/16}}},
  \label{eq:oracle_complexity_better}
\end{align}
where the $\O$-notation in the second term ignores dependencies on parameters other than $d$ and $\epsilon$.

\begin{algorithm}[t]
  \caption{Doubling scheme for \cref{alg:proposed}}
  \label{alg:doubling}
  \begin{algorithmic}[1]
    \Require{%
      $\xinit \in \R^d$,
      $c_\kappa > d^{1/5}$, $c_\sigma, c_\delta > 0$
    }
    \For{$T = 2^0, 2^1, 2^2, \dots$}
      \State Run \cref{alg:proposed} for $T$ outer iterations using the parameters in \cref{eq:parameter_schedule_given_T}
      \State Update $x_{\mathrm{init}}$ to the observed point with the lowest objective value
    \EndFor
  \end{algorithmic}
\end{algorithm}

The above analysis requires that the number of outer iterations $T$ is known in advance.
This requirement can be removed by using a doubling scheme, as described in \cref{alg:doubling}.

\begin{corollary}
  The oracle complexity of \cref{alg:doubling} is also given by \cref{eq:oracle_complexity_given_T}.
\end{corollary}
\begin{proof}
  Let $T_\epsilon$ be the number of outer iterations sufficient for \cref{alg:proposed} to find an $\epsilon$-stationary point.
  Let $m$ be the smallest integer such that $2^m \geq T_\epsilon$.
  Since the initial point is updated to the observed point with the lowest objective value, the value of $\Delta$ in each phase is at most the initial one.
  Hence, \cref{thm:oracle_complexity_given_T} applies to each $T = 2^0,\dots,2^m$.
  The total number of inner iterations is bounded by
  \begin{align}
    \sum_{i=0}^m c_\kappa \prn*{2^i}^{\frac{13}{12}}
    =
    c_\kappa \frac{2^{\frac{13}{12}(m+1)} - 1}{2^{\frac{13}{12}} - 1}
    =
    \O \prn*{c_\kappa 2^{\frac{13}{12} m}}
    =
    \O \prn*{c_\kappa T_\epsilon^{\frac{13}{12}}}.
  \end{align}
  Combining this bound with \cref{eq:num_outer_iterations} gives \cref{eq:oracle_complexity_given_T}.
\end{proof}

\section{Numerical experiments}
This section presents numerical experiments to evaluate the performance of \cref{alg:doubling}.
The parameter settings for the algorithm are described in \cref{sec:parameter_settings_experiments}.
All experiments were conducted using Python~3.13.9 on a MacBook Air equipped with an Apple M3 chip and 24~GB of memory.
The source code is available at \url{https://github.com/n-marumo/pf-aqnewton}.

\subsection{Benchmark functions}
We first investigate the performance on the following benchmark functions:
\begin{align}
  \text{Dixon--Price \citep{dixon1989truncated}:}\ \ 
  f(x)
  &= (x_1 - 1)^2 + \sum_{i=2}^d i (2 x_i^2 - x_{i-1})^2,
  \label{eq:dixon_price}
  \\
  \label{eq:powell}
  \text{Powell \citep{powell1962iterative}:}\ \ 
  f(x)
  &= \sum_{i=1}^{\floor{d/4}} \!\Big(
    \prn*{x_{4i-3} + 10 x_{4i-2}}^2
    + 5 \prn*{x_{4i-1} - x_{4i}}^2\\
    &\qquad\qquad
    + \prn*{x_{4i-2} - 2 x_{4i-1}}^4
    + 10 \prn*{x_{4i-3} - x_{4i}}^4
  \Big),
  \\
  \text{Qing \citep{qing2006dynamic}:}\ \
  f(x)
  &= \sum_{i=1}^d (x_i^2 - i)^2,
  \label{eq:qing}
  \\
  \text{Rosenbrock \citep{rosenbrock1960automatic}:}\ \
  f(x)
  &= \sum_{i=1}^{d-1} \prn*{
    100 \prn*{x_{i+1} - x_i^2}^2
    + (x_i - 1)^2
  },
  \label{eq:rosenbrock}
\end{align}
where $x = (x_1, \ldots, x_d) \in \R^d$.
The proposed method is evaluated against the following baseline algorithms.
\begin{itemize}
  \item 
  \BFGS and \DFP: standard quasi-Newton methods with Armijo's backtracking line search.
  We employ the cautious update strategy~\citep{li2001global}, which updates the Hessian approximation only when the curvature condition $\inner*{\nabla f(x_{k+1}) - \nabla f(x_k)}{x_{k+1} - x_k} > 0$ holds; otherwise the update is skipped to preserve positive definiteness.
  \item
  \GD: the gradient descent method with Armijo's backtracking line search.
  This method has an oracle complexity bound of $\O(\epsilon^{-2})$.
  \item
  \AGDR: the accelerated gradient descent with restart schemes~\citep{marumo2024parameter}.
  The parameter settings follow \citep{marumo2024parameter}.
  This method has an oracle complexity bound of $\O(\epsilon^{-7/4})$.
\end{itemize}

\Cref{table:experiments__benchmark_algos} summarizes the results on the benchmark functions with dimension $d = 100$.
For each function, we draw ten initial points independently from $\mathcal{N}(x^*, 10^4 I)$, where $x^*$ is the global minimizer, and run every algorithm from the same ten points.
The table reports the average number of oracle calls required to reach $\norm*{\nabla f(x)} \leq \epsilon$.
The results show that the proposed method is competitive with the baseline algorithms on the Dixon--Price, Qing, and Rosenbrock functions, particularly for smaller values of $\epsilon$.
On the Powell function, however, \BFGS consistently requires fewer oracle calls than the proposed method.
This indicates an empirical limitation of the proposed method on this problem.

\begin{table}[t]
  \centering
  \small
  \caption{%
    Average number of oracle calls required to reach $\norm*{\nabla f(x)} \leq \epsilon$.
    An entry marked with ``--'' indicates that the corresponding accuracy was not reached within $10^5$ oracle calls in at least one run.
    \label{table:experiments__benchmark_algos}
  }
  \begin{tabular}{@{}llrrrrr@{}}
\toprule
Function & Method & $\epsilon = 10^{1}$ & $10^{0}$ & $10^{-1}$ & $10^{-2}$ & $10^{-3}$ \\
\midrule
Dixon--Price & \Proposed & \textbf{923.0} & \textbf{1031.2} & \textbf{1073.2} & \textbf{1142.9} & \textbf{1283.6} \\
 & \BFGS & 1312.2 & 1661.2 & 1908.2 & 2071.2 & 2195.8 \\
 & \DFP & -- & -- & -- & -- & -- \\
 & \GD & 1073.4 & 1222.4 & 1491.5 & 1911.8 & 2445.2 \\
 & \AGDR & 15960.4 & 19805.8 & 22576.2 & 24421.8 & 26206.8 \\
\midrule
Powell & \Proposed & 2406.6 & 2488.6 & 2544.6 & 2622.4 & 2716.4 \\
 & \BFGS & \textbf{729.7} & \textbf{874.6} & \textbf{1037.0} & \textbf{1208.3} & \textbf{1353.1} \\
 & \DFP & 24431.4 & -- & -- & -- & -- \\
 & \GD & 806.5 & 1136.8 & 2320.6 & 7541.5 & 31228.9 \\
 & \AGDR & 9666.9 & 12041.9 & 15438.9 & 21132.1 & 28783.7 \\
\midrule
Qing & \Proposed & 198.8 & 391.5 & \textbf{413.7} & \textbf{443.3} & \textbf{509.9} \\
 & \BFGS & 666.2 & 829.4 & 930.0 & 998.6 & 1056.4 \\
 & \DFP & 5404.1 & 11223.5 & 12904.1 & 14271.9 & 14932.7 \\
 & \GD & \textbf{129.9} & \textbf{297.6} & 563.1 & 866.1 & 1178.1 \\
 & \AGDR & 777.7 & 1076.3 & 1354.3 & 1609.5 & 1988.7 \\
\midrule
Rosenbrock & \Proposed & \textbf{424.0} & \textbf{1540.4} & \textbf{1906.2} & \textbf{1968.8} & \textbf{2023.5} \\
 & \BFGS & 2131.8 & 3819.6 & 4046.8 & 4051.4 & 4055.2 \\
 & \DFP & -- & -- & -- & -- & -- \\
 & \GD & 6008.9 & 46560.4 & 53743.3 & 64568.8 & 76007.5 \\
 & \AGDR & 6405.3 & 27517.3 & 28856.3 & 30277.9 & 31145.8 \\
\bottomrule
\end{tabular}

\end{table}

\begin{figure}[t]
  \centering
  \includegraphics[height=3.2ex]{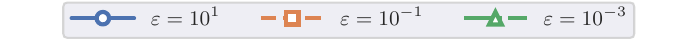}\par\medskip%
  \subfloat[Dixon--Price\label{fig:exp_dim_oracle_dixonprice}]{%
    \includegraphics[width=0.4\linewidth]{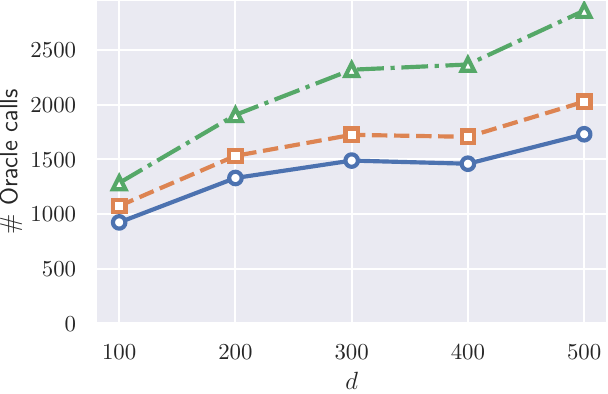}%
  }\hspace{0.1\linewidth}%
  \subfloat[Powell\label{fig:exp_dim_oracle_powell}]{%
    \includegraphics[width=0.4\linewidth]{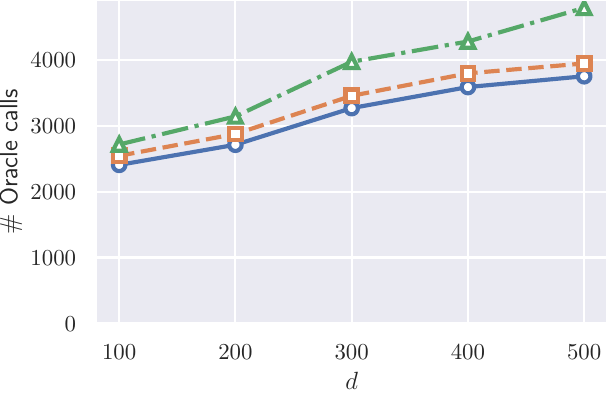}%
  }\par\smallskip%
  \subfloat[Qing\label{fig:exp_dim_oracle_qing}]{%
    \includegraphics[width=0.4\linewidth]{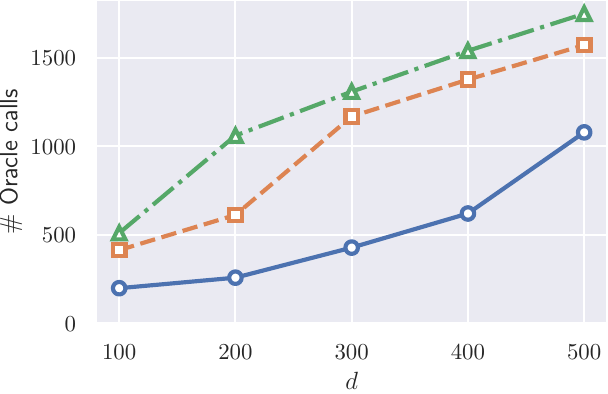}%
  }\hspace{0.1\linewidth}%
  \subfloat[Rosenbrock\label{fig:exp_dim_oracle_rosenbrock}]{%
    \includegraphics[width=0.4\linewidth]{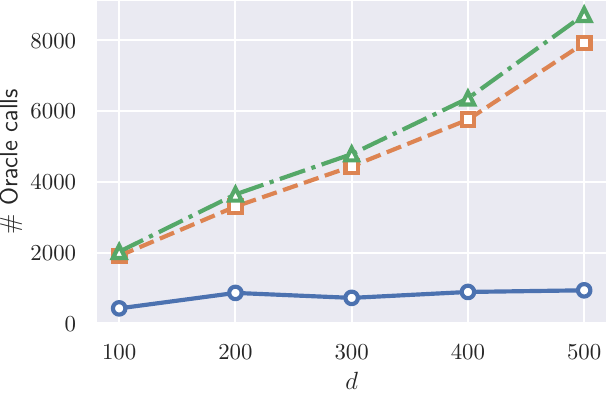}%
  }\par
  \caption{
    Average number of oracle calls required by the proposed method to reach $\norm*{\nabla f(x)} \leq \epsilon$.
    \label{fig:experiments_dim_oracle}
  }
\end{figure}

We also investigate how the proposed method scales with the problem dimension $d$ by testing it on benchmark functions with dimensions $d \in \set*{100, 200, 300, 400, 500}$.
For each dimension, we use ten random initial points generated in the same manner as those used for \cref{table:experiments__benchmark_algos}, and \cref{fig:experiments_dim_oracle} plots the average number of oracle calls.
Although the worst-case oracle complexity is $\O(d^{1/4} \epsilon^{-13/8})$, the plots should not be interpreted as evidence of a $\Theta(d^{1/4})$ dependence on $d$.
The benchmark functions are not necessarily worst-case instances, and the constants hidden in the $\O$-notation depend on $L$, $M$, and $\Delta$, which may also vary with $d$.

\subsection{Inverse problem for the wave equation}
We next consider an application-oriented instance, an inverse problem for the following nonlinear wave equation from \citep{bellavia2018levenberg,marumo2025accelerated}:
\begin{align}
  \frac{\partial^2 u}{\partial t^2} (z, t)
  =
  \frac{\partial^2 u}{\partial z^2} (z, t)
  - \e^{u(z, t)}
  \quad
  \text{for all $(z, t) \in [0, 1]^2$.}
\end{align}
The function $u$ satisfies $u(0, t) = u(1, t) = 0$ for all $t \in [0, 1]$ and $\frac{\partial u}{\partial t} (z, 0) = 0$ for all $z \in [0, 1]$.
The goal is to reconstruct the initial state $u(z, 0)$ from the final state $u(z, 1)$.
To formulate this problem, we discretize the wave equation using central differences (see, e.g., \citep[Chapter~6]{strang2007computational}) with grid points $(z, t) = (i/N, j/T)$ for $0 \leq i \leq N$ and $0 \leq j \leq T$.
For $w \in \R^{N-1}$, let $\hat u_i(w)$ be the approximation of $u(i/N, 1)$ obtained by solving the discretized equation with the initial state $u(i/N, 0) = w_i$ for $1 \leq i \leq N-1$.
This leads to the following optimization problem:
\begin{align}
  \min_{w \in \R^{N-1}}
  \frac{1}{N - 1} \sum_{i=1}^{N-1}
  \prn[\big]{ \hat u_i(w) - \hat u_i(w^\star) }^2,
\end{align}
where $w^\star \in \R^{N-1}$ denotes the discretized true initial state.
We set $N = 100$ and $T = 1000$, and use the same true initial state as in \citep{marumo2025accelerated}:
\begin{align}
  u^\star(z, 0)
  =
  \sin(6 \pi z)
  +
  \begin{dcases*}
    1 - \cos(20 \pi z) & if $\frac{4}{10} \leq z \leq \frac{5}{10}$, \\
    0 & otherwise.
  \end{dcases*}
\end{align}
We then set $w^\star = (u^\star(i/N, 0))_{1 \leq i \leq N-1}$.
The initial point is set to $w = \0$.

\begin{figure}[t]
  \centering
  \includegraphics[height=3.2ex]{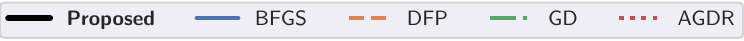}\par\medskip%
  \includegraphics[width=0.48\linewidth]{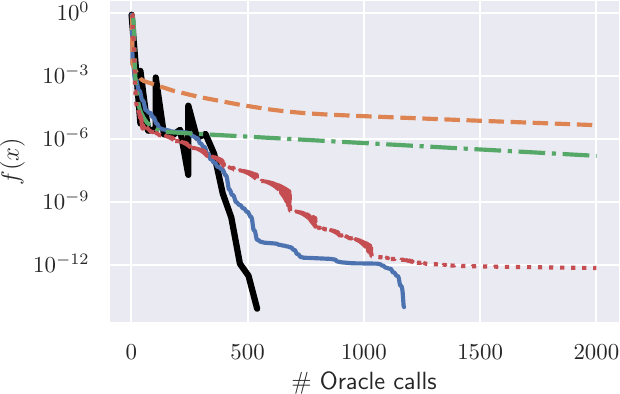}%
  \hspace{0.03\linewidth}%
  \includegraphics[width=0.48\linewidth]{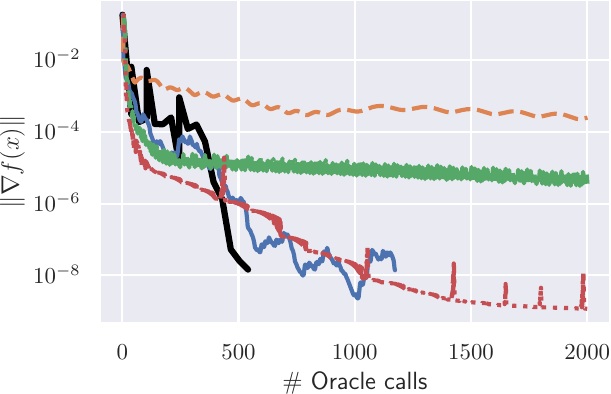}\par
  \caption{
    Results of the inverse problem for the nonlinear wave equation.
    \label{fig:experiments_wave}
  }
\end{figure}

\Cref{fig:experiments_wave} shows the results.
The proposed method is competitive with the baseline algorithms.
The results suggest that it is particularly promising when high-accuracy solutions are required.

\section{Discussion}
We proposed a quasi-Newton-type method that achieves an oracle complexity of $\O(d^{1/4} \epsilon^{-13/8})$.
The method is parameter-free in the sense that it does not require prior knowledge of problem-dependent parameters such as Lipschitz constants.
Moreover, it does not need to know the target accuracy $\epsilon$ or the total number of iterations in advance.

We now discuss how the complexity bounds depend on parameters other than $\epsilon$ and $d$.
If the algorithm parameters are chosen using knowledge of $L$, $M$, and $\Delta$, as in \cref{eq:setting_cs_with_d_L_M_Delta}, then the complexity bound~\cref{eq:oracle_complexity_better} is obtained.
On the other hand, under the setting \cref{eq:def_cs_only_d}, which is the main focus of this paper and does not use problem-dependent parameters, the dependence on $L$, $M$, and $\Delta$ becomes worse, as shown in \cref{eq:oracle_complexity_only_d}.
Such a tradeoff is also observed in other problem settings.
For example, in nonconvex stochastic optimization, an SGD-type algorithm can achieve a complexity of $\O(\Delta L \sigma^2 \epsilon^{-4})$ when tuned using problem-dependent parameters, whereas parameter-free tuning worsens the bound to $\O((\Delta^4 + L^4 + \sigma^4) \epsilon^{-4})$~\citep{hubler2025gradient}, where $\sigma^2$ denotes the variance of the stochastic gradient.
Improving the dependence on $L$, $M$, and $\Delta$ in our parameter-free complexity bound~\cref{eq:oracle_complexity_only_d} is left for future work.

We next discuss the computational complexity.
As discussed in \cref{sec:solving_subproblem}, solving each subproblem incurs a computational cost of $\tilde \O(d^2)$.
This cost is not an issue when $d$ is small or when a gradient evaluation is more expensive than solving the subproblem.
For large-scale problems, however, reducing this computational cost is desirable.
For example, in standard BFGS, the cost of solving the subproblem is $\O(d^2)$, whereas its limited-memory variant reduces the cost to $\O(md)$, where $m \geq 1$ is the memory size.
Developing a limited-memory variant of the proposed algorithm would be an interesting direction for future work.

\appendix

\section{Deferred proofs}
\subsection{Proof of \texorpdfstring{\cref{lem:gammak_diff_upperbound_general}}{Lemma~\ref{lem:gammak_diff_upperbound_general}}}
\label{sec:proof_lem_gammak_diff_upperbound_general}

\begin{proof}
  To prove the desired inequality, we rewrite $\inner*{X}{G_k}$ using the definition of $B_{k+1}$ in \cref{eq:update_Bk_general}:
  \begin{align}
    \inner*{X}{G_k}
    &=
    \frac{1+\theta}{1-\theta}\inner*{B_{k+1}}{G_k}
    - \inner*{B_k}{G_k},
    \label{eq:decomposition_inner_XG}
  \end{align}
  and then bound each term on the right-hand side.
  The first term is bounded as follows:
  \begin{align}
    \label{eq:upperbound_inner_BG1_pre}
    &\mathInd
    2 \inner*{B_{k+1}}{G_k}\\
    &=
    \frac{1}{1 - \theta} \normFro*{B_{k+1}}^2
    - \gamma_{k+1}
    + 2 \inner*{B_{k+1}}{G_k - \nabla^2 f(x_{k+1})}\\
    &\leq
    \frac{1}{1 - \theta} \normFro*{B_{k+1}}^2
    - \gamma_{k+1}
    + \frac{\theta}{1 - \theta} \normFro*{B_{k+1}}^2
    + \frac{1 - \theta}{\theta} \normFro*{G_k - \nabla^2 f(x_{k+1})}^2\\
    &=
    \frac{1 + \theta}{1 - \theta} \normFro*{B_{k+1}}^2
    - \gamma_{k+1}
    + \frac{1 - \theta}{\theta} \normFro*{G_k - \nabla^2 f(x_{k+1})}^2,
  \end{align}
  where the first equality follows from the equality in \cref{eq:gamma_equality_lowerbound} and the inequality follows from Young's inequality.
  The last term in \cref{eq:upperbound_inner_BG1_pre} can be further bounded as
  \begin{align}
    \normFro*{G_k - \nabla^2 f(x_{k+1})}
    &\leq
    \sqrt{d} \normop*{G_k - \nabla^2 f(x_{k+1})}\\
    &\leq
    \sqrt{d} \int_0^1 \normop*{\nabla^2 f(x_k + \tau s_k) - \nabla^2 f(x_{k+1})} \,\dd \tau\\
    &\leq
    \sqrt{d} \int_0^1 M (1 - \tau) \norm*{s_k} \,\dd \tau
    =
    \frac{\sqrt{d} M}{2} \norm*{s_k},
  \end{align}
  where the second inequality follows from the definition of $G_k$ in \cref{eq:def_Gk} and the last inequality follows from $x_{k+1} = x_k + s_k$ and \cref{asm:lipschitz_hess}.
  Thus, we have
  \begin{align}
    2 \inner*{B_{k+1}}{G_k}
    &\leq
    \frac{1 + \theta}{1 - \theta} \normFro*{B_{k+1}}^2
    - \gamma_{k+1}
    + \frac{1 - \theta}{4 \theta} d M^2 \norm*{s_k}^2.
    \label{eq:upperbound_inner_BG1}
  \end{align}
  \FloatBarrier
  The second term on the right-hand side of \cref{eq:decomposition_inner_XG} is bounded in the same manner:
  \begin{align}
    \label{eq:upperbound_inner_BG2}
    \quad
    - 2 \inner*{B_k}{G_k}
    &=
    \gamma_k
    - \frac{1}{1 - \theta} \normFro*{B_k}^2
    + 2 \inner*{B_k}{\nabla^2 f(x_k) - G_k}\\
    &\leq
    \gamma_k
    - \frac{1}{1 - \theta} \normFro*{B_k}^2
    + \frac{\theta}{1 - \theta} \normFro*{B_k}^2
    + \frac{1 - \theta}{\theta} \normFro*{\nabla^2 f(x_k) - G_k}^2\\
    &\leq
    \gamma_k
    - \normFro*{B_k}^2
    + \frac{1 - \theta}{4 \theta} d M^2 \norm*{s_k}^2.
  \end{align}
  Plugging \cref{eq:upperbound_inner_BG1,eq:upperbound_inner_BG2} into \cref{eq:decomposition_inner_XG} yields
  \begin{align}
    2 \inner*{X}{G_k}
    &\leq
    \frac{1 + \theta}{1 - \theta} \prn*{
      \frac{1 + \theta}{1 - \theta} \normFro*{B_{k+1}}^2
      - \gamma_{k+1}
      + \frac{1 - \theta}{4 \theta} d M^2 \norm*{s_k}^2
    }\\
    &\qquad
    + \prn*{
      \gamma_k
      - \normFro*{B_k}^2
      + \frac{1 - \theta}{4 \theta} d M^2 \norm*{s_k}^2
    }\\
    &=
    \gamma_k - \gamma_{k+1}
    + \frac{(1 + \theta)^2}{(1 - \theta)^2} \normFro*{B_{k+1}}^2
    - \normFro*{B_k}^2
    + \frac{d M^2}{2 \theta} \norm*{s_k}^2
    - \frac{2 \theta}{1 - \theta} \gamma_{k+1}. 
  \end{align}
  
  Taking the Frobenius norm of both sides of \cref{eq:update_Bk_general} gives
  \begin{align}
    \frac{(1 + \theta)^2}{(1 - \theta)^2} \normFro*{B_{k+1}}^2
    &=
    \normFro*{B_k}^2
    + \normFro*{X}^2
    + 2 \inner*{X}{B_k}.
  \end{align}
  Substituting this relation into the previous inequality and rearranging terms, we obtain
  \begin{align}
    \gamma_{k+1} - \gamma_k
    &\leq 
    \frac{d M^2}{2 \theta} \norm*{s_k}^2
    + \normFro*{X}^2
    - 2 \inner*{X}{G_k - B_k}
    - \frac{2 \theta}{1 - \theta} \gamma_{k+1}.
  \end{align}
  Finally, applying \cref{eq:gamma_equality_lowerbound} to the last term gives $- \frac{2 \theta}{1 - \theta} \gamma_{k+1} \leq  2 \theta d L^2$, which completes the proof.
\end{proof}

\subsection{Proof of \texorpdfstring{\cref{lem:optimal_X}}{Lemma~\ref{lem:optimal_X}}}
\label{sec:proof_lem_optimal_X}
\begin{proof}
  Since $G_k s_k = \nabla f(x_{k+1}) - \nabla f(x_k)$, note that
  \[
    \prn*{G_k - B_k} s_k
    =
    \nabla f(x_{k+1}) - \nabla f(x_k) - B_k s_k
    =
    r_k.
  \]
  Rewriting the problem~\cref{eq:problem_for_X} gives
  \begin{align}
    \label{eq:reformulated_problem_for_X}
    &\mathInd
    {\min_{v \in \R^d}} \
    \set*{
      \big\| v s_k^\top + s_k v^\top \big\|_{\mathrm{F}}^2
      - 2 \inner*{v s_k^\top + s_k v^\top}{G_k - B_k}
    }\\
    &=
    \min_{v \in \R^d} \
    \set*{
      2 \norm*{s_k}^2 \norm*{v}^2
      + 2 \inner*{v}{s_k}^2
      - 4 \inner*{v}{r_k}
    }.
  \end{align}
  The first-order optimality condition yields
  \[
    4 \norm*{s_k}^2 v
    + 4 \inner*{v}{s_k} s_k
    - 4 r_k
    = 0.
  \]
  Solving this equation for $v$ gives the unique optimal solution:
  \[
    v
    =
    \prn*{\norm*{s_k}^2 I + s_k s_k^\top}^{-1} r_k
    =
    \prn*{
      \frac{1}{\norm*{s_k}^2}I - \frac{1}{2 \norm*{s_k}^4}s_k s_k^\top
    } r_k
    =
    \frac{r_k}{\norm*{s_k}^2} - \frac{\inner*{r_k}{s_k}}{2 \norm*{s_k}^4} s_k,
  \]
  where the second equality follows from the Sherman--Morrison formula.
  Substituting this optimal $v$ into $X = v s_k^\top + s_k v^\top$ yields the optimal $X$ in \cref{eq:optimal_X}.

  Substituting the optimal $v$ into the objective function in \cref{eq:reformulated_problem_for_X}, we obtain the optimal value as follows:
  \begin{align}
    &\mathInd
    2 \norm*{s_k}^2 \norm*{\frac{r_k}{\norm*{s_k}^2} - \frac{\inner*{r_k}{s_k}}{2 \norm*{s_k}^4} s_k}^2
    + 2 \inner*{\frac{r_k}{\norm*{s_k}^2} - \frac{\inner*{r_k}{s_k}}{2 \norm*{s_k}^4} s_k}{s_k}^2\\
    &\qquad
    - 4 \inner*{\frac{r_k}{\norm*{s_k}^2} - \frac{\inner*{r_k}{s_k}}{2 \norm*{s_k}^4} s_k}{r_k}\\
    &=
    2 \norm*{s_k}^2
    \prn*{
      \frac{\norm*{r_k}^2}{\norm*{s_k}^4}
      - \frac{3}{4} \frac{\inner*{r_k}{s_k}^2}{\norm*{s_k}^6}
    }
    + 2 \prn*{ \frac{\inner*{r_k}{s_k}}{2 \norm*{s_k}^2} }^2
    - 4 \prn*{
      \frac{\norm*{r_k}^2}{\norm*{s_k}^2}
      - \frac{\inner*{r_k}{s_k}^2}{2 \norm*{s_k}^4}
    }\\
    &=
    \frac{\inner*{r_k}{s_k}^2}{\norm*{s_k}^4}
    - 2 \frac{\norm*{r_k}^2}{\norm*{s_k}^2},
  \end{align}
  which completes the proof.
\end{proof}

\subsection{Proof of \texorpdfstring{\cref{lem:grad_norm_at_average}}{Lemma~\ref{lem:grad_norm_at_average}}}
\label{sec:proof_lem_grad_norm_at_average}
\begin{proof}
  From \citep[Lemma~3.1]{marumo2024parameter}, for any $\alpha_0,\dots,\alpha_k \geq 0$ with $\sum_{i=0}^k \alpha_i = 1$, we have
  \[
    \norm*{\nabla f\prn*{ \sum_{i=0}^k \alpha_i x_i } - \sum_{i=0}^k \alpha_i \nabla f(x_i)}
    \leq
    \frac{M}{2} \sum_{0 \leq i < j \leq k} \alpha_i \alpha_j \norm*{x_i - x_j}^2.
  \]
  Set $\alpha_i = \frac{2i+1}{k (k+1)}$ for $0 \leq i < k$ and $\alpha_k = \frac{1}{k+1}$.
  Using the definitions of $\bar x_k$ and $\bar g_k$ in \cref{eq:def_xbar,eq:def_gbar} yields
  \begin{align}
    \norm*{\nabla f\prn*{ \bar x_k } - \bar g_k}
    &\leq
    \frac{M}{2} \sum_{0 \leq i < j \leq k} \alpha_i \alpha_j \norm*{x_i - x_j}^2.
    \label{eq:grad_norm_proof1}
  \end{align}
  Furthermore, for $0 \leq i < j \leq k$, it follows from the triangle inequality and the Cauchy--Schwarz inequality that
  \[
    \norm*{x_i - x_j}^2
    \leq
    \prn*{ \sum_{l=i}^{j-1} \norm*{s_l} }^2
    \leq
    \prn*{ \sum_{l=i}^{j-1} 1^2 }
    \prn*{ \sum_{l=i}^{j-1} \norm*{s_l}^2 }
    \leq
    (j - i) S_k.
  \]
  Plugging this into \cref{eq:grad_norm_proof1} gives
  \begin{align}
    \norm*{\nabla f(\bar x_k) - \bar g_k}
    &\leq
    \frac{M S_k}{2 k^2 (k+1)^2}
    \prn[\Bigg]{
      \sum_{0 \leq i < j < k} \!\!\! (2i+1)(2j+1) (j - i)
      + k \sum_{i=0}^{k-1} (2i+1) (k - i)
    }\\
    &=
    \frac{M S_k}{2 k^2 (k+1)^2}
    \frac{k(k+1)(2k+1)(2k^2 + 2k + 1)}{30}\\
    &\leq
    \frac{M S_k}{2 k^2 (k+1)^2}
    \frac{k(k+1)(2k+k)(2k^2 + 2k + \frac{k^2 + k}{2})}{30}
    =
    \frac{M}{8} k S_k.
  \end{align}
  Applying the triangle inequality concludes the proof.
\end{proof}

\section{Approximate eigendecomposition of $B_k$}
\label{sec:approx_eigendecomposition}
To apply \cref{lem:subproblem_reduction} to the subproblem~\cref{eq:subproblem}, we maintain the approximation condition~\cref{eq:approx_eigendecomposition_B} for each $B = B_k$.
Specifically, we store the exact matrix $B_k$ together with an eigendecomposition of its approximation $\hat B_k$:
\begin{align}
  \hat B_k
  \coloneqq
  Q_k \Lambda_k Q_k^\top,
  \quad
  Q_k^\top Q_k
  =
  I,
  \quad
  \normop[\big]{B_k - \hat B_k}
  \leq
  \frac{\delta}{2}.
  \label{eq:auxiliary_eigendecomposition_Bk}
\end{align}
The scaled PSB update~\cref{eq:update_Bk} consists of a rank-two update followed by scaling, and the rank-two update can be decomposed into two rank-one updates.
If each rank-one update introduces error at most $\frac{\theta \delta}{2 \prn*{1 - \theta}}$, then combining \cref{eq:update_Bk,eq:auxiliary_eigendecomposition_Bk} gives
\[
  \normop[\big]{B_{k+1} - \hat B_{k+1}}
  \leq
  \frac{1 - \theta}{1 + \theta} \prn*{
    \normop[\big]{B_k - \hat B_k}
    + \frac{\theta \delta}{1 - \theta}
  }
  \leq
  \frac{1 - \theta}{1 + \theta} \prn*{
    \frac{\delta}{2}
    + \frac{\theta \delta}{1 - \theta}
  }
  =
  \frac{\delta}{2}.
\]
Rank-one updates of eigendecompositions are classical and can be computed to fixed accuracy with $\tilde{\O} (d^2)$ computational cost~\citep{gu1994stable,gu1995divide}.
More specifically, the Bunch--Nielsen--Sorensen algorithm~\citep{bunch1978rank} computes the updated eigenvalues, and fast Cauchy matrix multiplication~\citep{gu1995divide} can then be used to compute the updated eigenvectors.

\section{Optimality of the exponents in the parameter settings}
\label{sec:parameter_optimality}
We examine the following parameter settings, which are more general than those in \cref{eq:parameter_schedule_given_T}:
\begin{align}
  \kappa_t
  =
  c_\kappa T^a,
  \qquad
  \sigma_t
  =
  c_\sigma T^b,
  \qquad
  \delta_t
  =
  c_\delta T^c,
  \label{eq:general_power_parameter_schedule}
\end{align}
where $c_\kappa > d^{1/5}$, $c_\sigma, c_\delta > 0$, $a > 0$, and $b, c \in \R$.
Plugging \cref{eq:general_power_parameter_schedule} into \cref{eq:gradnorm_bound_general} gives $\min_{0 \leq t < T} \norm{\nabla f \prn{\bar x_{K_t}^t}} \leq \O \prn*{T^{-q}}$, where
\begin{align}
  q
  \coloneqq
  \min
  \set*{
    1 + 4a - b,\,
    1 + 4a,\,
    8a,\,
    3a - 2c,\,
    2b - 8a
  }.
  \label{eq:def_general_power_rate}
\end{align}
Since the total number of inner iterations is $\O \prn*{T^{1+a}}$, the resulting oracle complexity is $\O \prn*{ \epsilon^{- (1+a)/q} }$.

We minimize the exponent $(1+a)/q$ to obtain the best dependence on $\epsilon$.
The definition of $q$ in \cref{eq:def_general_power_rate} implies $q \leq 1 + 4a - b$ and $q \leq 2b - 8a$.
Multiplying the first inequality by $2$ and adding it to the second inequality gives $3q \leq 2$.
The definition of $q$ also implies $q \leq 8a$.
Combining these bounds yields
\begin{align}
  \frac{1 + a}{q}
  = \frac{1}{q} + \frac{a}{q}
  \geq
  \frac{3}{2} + \frac{1}{8}
  =
  \frac{13}{8}.
\end{align}
The lower bound is attained if and only if
\begin{align}
  a
  =
  \frac{1}{12},
  \quad
  b
  =
  \frac{2}{3},
  \quad
  c
  \leq
  - \frac{5}{24}.
\end{align}
Therefore, the exponents of $\kappa_t$ and $\sigma_t$ in \cref{eq:parameter_schedule_given_T} are optimal within the choices in \cref{eq:general_power_parameter_schedule}.
Any exponent $c \leq -5/24$ for $\delta_t$ retains the complexity $\O \prn*{\epsilon^{-13/8}}$, while $c = -5/24$ permits the largest tolerance for subproblem solutions.

\section{Parameter settings for numerical experiments}
\label{sec:parameter_settings_experiments}

\Cref{alg:doubling} has three constants: $c_\kappa$, $c_\sigma$, and $c_\delta$.
We fix $c_\delta = 10^{-5}$ because $c_\delta$ has negligible influence on the algorithm as long as it is sufficiently small.

We next discuss the choice of $c_\kappa$, whose role is easy to interpret: it controls the number of inner iterations $K$.
Indeed, from Line~\ref{alg-line:set_parameters} of \cref{alg:proposed} and \cref{eq:parameter_schedule_given_T}, we have $K = \floor*{c_\kappa T^{1/12}}$.
Since the exponent $1/12$ is small, $K$ grows very slowly with $T$, and its practical scale is mainly determined by $c_\kappa$.
Guided by the dimension dependence in \cref{eq:def_cs_only_d}, we set $c_\kappa = 10 d^{1/4}$ in the numerical experiments.
The constant factor $10$ is used to avoid an excessively small number of inner iterations.
For $d = 100$ and $1 \leq T \leq 2^{12}$, this choice yields $31 \leq K \leq 63$.

Finally, we describe the choice of $c_\sigma$, which controls the regularization parameter $\sigma$.
Since $\sigma$ should be chosen according to the scale of the problem, the choice of $c_\sigma$ requires more care.
We set $c_\sigma$ by the following procedure, inspired by the theoretical choice in \cref{eq:setting_cs_with_d_L_M_Delta}.
First, to obtain an estimate $\hat L$ of $L$, we perform a backtracking line search for a gradient-descent step at the initial point $x_{\mathrm{init}}$.
Specifically, with $x = x_{\mathrm{init}}$ and $y = x - \frac{1}{\hat L} \nabla f(x)$, we find $\hat L > 0$ by backtracking so that
\begin{align}
  f(y) - f(x)
  \leq
  \inner*{\nabla f(x)}{y - x}
  + \frac{\hat L}{2} \norm*{y - x}^2
  =
  - \frac{1}{2 \hat L} \norm*{\nabla f(x)}^2.
\end{align}
Once $y$ is determined, we estimate $M$ using a Hessian-free inequality of the same form as \cref{eq:hessian-free-ineq-function}:
\begin{align}
  \hat M
  \coloneqq
  \frac{12}{\norm{y - x}^3} \abs*{ f(y) - f(x) - \frac{1}{2} \inner*{\nabla f(x) + \nabla f(y)}{y - x} }.
\end{align}
We also estimate $\Delta$ by $\hat \Delta \coloneqq f(x) - f(y)$.
Using these estimates, and following \cref{eq:setting_cs_with_d_L_M_Delta}, we set $c_\sigma = \frac{d \hat L \hat M^{2/3}}{\hat \Delta^{2/3}}$.

\section*{Acknowledgments}
The author is grateful to Akiko Takeda for her valuable comments on the manuscript and to the anonymous reviewers for their constructive comments and suggestions.

\bibliographystyle{siamplain}
\bibliography{ref_siopt}
\end{document}